\documentclass[11pt,a4paper]{article}
\usepackage[utf8]{inputenc}
\usepackage{microtype}
\usepackage{amsmath,amsthm,amssymb}
\usepackage{mathrsfs} 
\usepackage[T1]{fontenc}
\usepackage{ae,aecompl}
\usepackage{slashed}
\usepackage{times}
\usepackage[british]{babel}
\usepackage[colorlinks,pagebackref,hypertexnames=false]{hyperref} 
\usepackage{url}
\usepackage[numeric,backrefs]{amsrefs}
\usepackage{nth}
\usepackage{float}
\usepackage[all]{xy}
\usepackage{tikz-cd} 
\usepackage{bookmark}
\usepackage[margin=2cm]{geometry}
\usepackage{booktabs}
\usepackage{tabularx}
\usepackage{enumerate}

\usepackage[textwidth=1.7cm,textsize=scriptsize]{todonotes}
\usepackage{comment}

\newtheorem{theorem}{Theorem}
\newtheorem{proposition}[theorem]{Proposition}

\newtheorem{lemma}[theorem]{Lemma}

\theoremstyle{remark} 
\newtheorem{remark}{Remark}

\theoremstyle{definition}
\newtheorem{definition}[theorem]{Definition}
\newtheorem{example}[theorem]{Example}

\numberwithin{theorem}{section}
\numberwithin{remark}{section}
\newcommand{\R}{\mathbb{R}}
\newcommand{\G}{\mathrm{G}}
\newcommand{\SO}{\mathrm{SO}}
\newcommand{\SU}{\mathrm{SU}}

\newcommand{\vp}{\varphi}
\newcommand{\w}{\wedge}

\newcommand{\Ooct}{\mathbb{O}}

\newcommand{\dd}{\mathrm{d}}
\newcommand{\Hol}{\operatorname{Hol}}
\newcommand{\Ric}{\operatorname{Ric}}

\newcommand{\ad}{\operatorname{ad}}
\newcommand{\vol}{\operatorname{vol}}

\allowdisplaybreaks
\date{}

\usepackage{titling}
\usepackage{changepage}

\usepackage{authblk}
\usepackage{enumitem}
\usepackage{mathtools}

\title{Strong $\G_2$-structures with torsion: a survey and new examples}

\author{Anna Fino,\qquad \qquad  Udhav Fowdar}
\date{}

\begin{document}
\maketitle

\begin{abstract}
This article has two main parts. The first consists of a brief survey of $\G_2$-structures with torsion, namely $\G_2$-structures admitting a $\G_2$-connection with totally skew-symmetric torsion. We focus primarily on the strong condition that the associated torsion $3$-form is closed. We review their basic properties and their connections with the heterotic $\G_2$ system, generalized geometry, and six-dimensional geometry. In the second part, we construct two complete cohomogeneity-one strong $\G_2T$-structures on $\mathbb{R}\times S^3\times S^3$ with non-zero closed and co-closed characteristic torsion, and non-vanishing characteristic Ricci tensor. Both have an asymptotically cylindrical end; at the other end, the first is asymptotically conical (AC), with a torsion-free $\G_2$ cone as its asymptotic model, while the second is asymptotically locally conical (ALC), with a constant-length circle bundle over a Calabi--Yau cone as its asymptotic model. We investigate the induced $\SU(3)$-structures on their circle quotients and show that they are globally conformally almost K\"ahler and globally conformally K\"ahler, respectively. To the best of our knowledge, these provide the first complete strong $\G_2T$-structures whose underlying metrics are not locally Riemannian products.

\end{abstract}
\begin{adjustwidth}{0.95cm}{0.95cm}
    \tableofcontents
\end{adjustwidth}

\addtocontents{toc}{\protect\setcounter{tocdepth}{1}} 
\section{Introduction}

Among the five exceptional simple Lie groups, the compact Lie group $\G_2$ plays a distinguished role in differential geometry. It occurs in Berger's list of possible irreducible Riemannian holonomy groups \cite{Berger} and is closely related to special holonomy, calibrated geometry, spin geometry and mathematical physics. A second, and particularly useful, point of view is to regard $\G_2$ as the stabiliser of a stable $3$-form on $\R^7$.

A {\em $\G_2T$-structure} on a $7$-manifold $M$ is a $\G_2$-structure $\varphi$ for
which there exists a metric connection $\nabla^T$ preserving $\varphi$ and
having totally skew-symmetric torsion $T_\varphi$.  When it exists, this
connection is unique and is called the characteristic connection.  The
$\G_2T$-structure is said to be \emph{strong} if its characteristic torsion is
closed, namely if $\dd T_\varphi=0$.  This is the seven-dimensional analogue of
the strong K\"ahler with torsion condition in Hermitian geometry.

This article has two complementary aims.  The first is to review the geometry
of $\G_2T$-structures, with particular emphasis on the strong case.  These
structures arise naturally in heterotic string compactifications, where the
torsion is interpreted as the NS-flux and is coupled to gauge fields through
the anomaly cancellation equation.  They also have a natural interpretation
in generalized geometry and are closely related to generalized Ricci
curvature.  We discuss rigidity phenomena, characteristic Ricci-flatness,
reduction to six-dimensional geometry and flows related to generalized Ricci
flow. The second aim is to construct two complete strong $\G_2T$-structures
with non-zero torsion and study their circle quotients.
The main  contribution of the paper is the following result.

\medskip
\noindent\textbf{Main result.}
\emph{There exist two complete strong $\G_2T$-structures on
$\R\times S^3\times S^3$ with non-zero harmonic characteristic
torsion. Both have an asymptotically cylindrical end.
The first is $\SU(2)^3$-invariant and has an asymptotically conical
end, while the second is $\SU(2)^2\times \mathrm{U}(1)$-invariant and has
an asymptotically locally conical end.  Both metrics have full Riemannian holonomy $\SO(7)$, and in both cases the characteristic Ricci tensor is non-zero.}\medskip

More precisely, we derive a four-function first-order system for
$\SU(2)^2\times \mathrm{U}(1)$-invariant strong $\G_2T$-structures and identify
two invariant two-dimensional subsystems containing complete
trajectories. At their cylindrical ends, these trajectories converge
to homogeneous product solutions studied in \cite{FinoFowdar2026};
at the other ends they have AC and ALC asymptotics, respectively.
Phase-plane arguments and ODE analysis give global existence, completeness and the asymptotic behaviour. The precise results are given in Theorem \ref{thm:mainsu(2)3} and \ref{thm:mainALC}.

In \cite{FinoFowdar2026}, we investigated the $\SU(2)^3$-invariant ansatz and searched for complete metrics which extend to the spinor bundle of $S^3$, namely, $S^3 \times \R^4$. It was shown that smooth completion at the singular $S^3$
orbit forces the torsion parameter to vanish leading to only torsion free examples. Our construction here instead achieves completeness by approaching a homogeneous
product solution at infinite distance backwards in the radial direction and yielding an AC  geometry in forward direction.
Passing to $\SU(2)^2\times \mathrm{U}(1)$ provides additional degrees of freedom and yields the ALC example.

{The circle quotients of these examples carry complete
$\SU(3)$-structures on $\R\times S^2\times S^3$.
They are globally conformally almost K\"ahler and globally
conformally K\"ahler, respectively. The latter has vanishing
Bismut-Ricci form, and its Bismut torsion satisfies
$\dd T_\omega=-F_\eta\wedge F_\eta,$
where $F_\eta$ is the curvature of the circle connection.}
The results and methods discussed here also point to several further
directions.  On the compact side, it would be useful to understand the global
structure of characteristic Ricci-flat strong $\G_2T$-manifolds beyond the
known product and homogeneous models.  The parallelism of the Lee form
provides a strong starting point, although the resulting six-dimensional
quotient geometry still has substantial freedom.

{For the cohomogeneity-one construction, a natural problem is to
determine the global phase portrait of the four-function system
derived in Section~7. More generally, one may ask whether complete
solutions with non-zero torsion and singular orbits can be obtained
within a larger $\SU(2)^2$-invariant ansatz.} From the analytic point of view, the
relation with generalized Ricci flow suggests looking for evolution equations
which preserve the strong condition and have characteristic Ricci-flat
structures as fixed points modulo diffeomorphisms.

{The paper is organised into two parts. Sections~2--6 survey
$\G_2$-structures with torsion and their relation to heterotic
string theory and generalized geometry. Sections~2 and~3 recall
intrinsic torsion and characteristic connections; Section~4
discusses the heterotic $\G_2$ system, Section~5 focuses on strong
$\G_2T$-structures, and Section~6 considers related geometric flows.
The second part contains the complete cohomogeneity-one
constructions in Section~7 and the study of their circle quotients
in Section~8. For more comprehensive introductions to exceptional
holonomy and $\G_2$ geometry, we refer to
\cites{Bryant1987,Bryant2006,JoyceI,JoyceII}.}

\smallskip
{\bf Acknowledgements.} 
The authors would like to thank Lorenzo Foscolo for useful discussions.
The first author is partially supported by  GNSAGA (INdAM) and by a grant from the Simons Foundation (\# 944448).

\smallskip

\textbf{AI disclosure.} This work was supported in part by a grant of access to OpenAI models through the ChatGPT for Academic Researchers program. While discussions with ChatGPT helped in refining some of the mathematical arguments, the article was itself written by the authors.

\section{\texorpdfstring{$\G_2$}{}-structures and intrinsic torsion}\label{sec: g2 structures and torsion}

There are several equivalent definitions of the Lie group $\G_2$. It is the automorphism
group of the octonion algebra $\Ooct$, it is the stabiliser of a positive $3$-form
on $\R^7$, and it is the stabiliser of a non-zero spinor in the real spin representation of $\operatorname{Spin}(7)$. These definitions lead to somewhat different, but nonetheless equivalent, approaches to $\G_2$-geometry: using octonions, differential forms, or spinors. We recall the basics and refer the reader to \cites{Baez2002, Bryant1987, Bonan1966, Joyce2000, Salamon1989} for more detail.

Identifying $\R^7\cong\operatorname{Im}\Ooct$, the octonionic multiplication defines the standard Euclidean metric $\langle \cdot,\cdot\rangle$ and a vector 
cross product $\times$ by
\[
 \langle x,y\rangle= -\mathrm{Re}(xy)\qquad\text{and} \qquad  x\times y =\mathrm{Im}(xy).
\]
This, in turn, defines a positive $3$-form $\varphi_0$ by
\begin{equation}\label{crossproduct}
 \varphi_0(x,y,z)=\langle x\times y,z\rangle.
\end{equation}
With respect to a suitable oriented orthonormal basis $\{e_i\}_{i=1}^7$, one may take
\begin{equation}\label{eq:model-form}
 \varphi_0=e^{127}+e^{347}+e^{567}+e^{135}-e^{245}-e^{146}-e^{236},
\end{equation}
where $e^{ijk}$ is shorthand for $e^i \wedge e^j\wedge e^k$. For instance, in this notation  we have $e_1 \times e_2 = e_7$. The $\mathrm{GL}(7,\R)$-orbit of $\varphi_0$ is open in $\Lambda^3(\R^7)^*$ and
\[
 \G_2 \cong  \{A\in \mathrm{GL}(7,\R):A^*\varphi_0=\varphi_0\}.
\]
\begin{definition}
A $\G_2$-structure on a $7$-manifold $M$ is a $3$-form $\varphi$ which is
pointwise equivalent to $\varphi_0$.
\end{definition}
The $3$-form $\varphi$ determines a Riemannian metric $g_\varphi$ and volume form $\vol_{\varphi}$ on $M$; and hence also a Hodge operator $*$. Following \cite{Bryant2006}, one convenient normalisation is
\[
 g_\varphi(X,Y)\vol_\varphi
 =\frac16\,\iota_X\varphi\wedge\iota_Y\varphi\wedge\varphi.
\]
In particular, this shows that $\G_2$ is a subgroup of $\SO(7)$. The equivalence of the first two descriptions now follows from \eqref{crossproduct} whereby $g_{\varphi}$ is identified with $\langle\cdot,\cdot\rangle$ at each point.

Since $\G_2$ is simply connected, it lifts to $\operatorname{Spin}(7)$ (the double cover of $\SO(7)$). As such, a $7$-manifold endowed with a $\G_2$-structure is spin, i.e. the first and second Stiefel–Whitney classes vanish. Moreover, $\varphi$ determines, up to sign, a unit real spinor $\eta$ as $S^7\cong \mathrm{Spin}(7)/\G_2$. Conversely, a unit spinor on a Riemannian $7$-manifold defines a $\G_2$-structure by
\begin{equation}\label{spinorial}
 \varphi(X,Y,Z)=\langle X\cdot Y\cdot Z\cdot\eta,\eta\rangle,
\end{equation}
where $\cdot$ denotes the Clifford multiplication. Note that the real spinor bundle is 8-dimensional and can be identified with $\Lambda^0\oplus\Lambda^1$ using the spinor $\eta$. The existence question is purely topological: a $7$-manifold admits a $\G_2$-structure if and only if it is spin \cites{Gray1969, LawsonMichelsohn1989}. Thus, \eqref{crossproduct} and \eqref{spinorial} exhibit the equivalence of the octonionic, differential-form, and spinorial descriptions of $\G_2$. 

Next we consider the differential geometric aspect of $\G_2$-structures. The torsion-free condition can be characterised in several equivalent ways.
\begin{theorem}[Fern\'andez--Gray \cite{FernandezGray1982}]
For a $\G_2$-structure defined by $\varphi$ the following are equivalent:
\begin{enumerate}[label=(\roman*)]
\item $\nabla^{LC}\varphi=0$, where $\nabla^{LC}$ denotes the Levi-Civita connection of the metric $g_{\varphi}$;
\item $\dd\varphi=0$ and $\dd *\varphi=0$;
\item $\Hol(g_\varphi)\subseteq \G_2$.
\end{enumerate}
\end{theorem}
If one, and hence any, of the above conditions hold then  the $\G_2$-structure is called \emph{torsion-free} or \emph{parallel}. In this case the induced metric $g_{\varphi}$ is automatically Ricci-flat \cite{Bonan1966}. In the language of spinors, the torsion free condition is equivalent to saying that the corresponding spinor $\eta$ is parallel with respect to the spin connection induced by the Levi-Civita connection.

The existence problem for torsion-free \(\G_2\)-structures is by now well developed. The first local examples with holonomy \textit{equal} to \(\G_2\) were constructed by Bryant \cite{Bryant1987}, complete non-compact examples by Bryant--Salamon \cite{BryantSalamon1989}, and the first compact examples by Joyce \cites{JoyceI,JoyceII}. We refer the reader to \cite{KarigiannisLeungLotay2020} for more recent developments in exceptional holonomy. 

In this article, however, we are primarily interested in \(\G_2\)-structures with non-zero intrinsic torsion.
Before giving a precise definition, we first recall some basic $\G_2$ representation theory. 

Under the action of $\G_2$, the spaces of $2$-forms and $3$-forms decompose into irreducible components as
\[
 \Lambda^2=\Lambda^2_7\oplus\Lambda^2_{14},\quad
 \Lambda^3=\Lambda^3_1\oplus\Lambda^3_7\oplus\Lambda^3_{27},
\]
where these irreducible components can be characterised as follows:
\[
\Lambda^2_7
=\{\iota_X \varphi,  \, X  \in  \Gamma(TM) \},
\quad
\Lambda^2_{14}
=\{\beta\in\Lambda^2:\beta\wedge *\varphi=0\},
\]
and
\[
\Lambda^3_1=\langle\varphi\rangle,
\quad
\Lambda^3_7
=\{\iota_X (* \varphi),  \, X  \in  \Gamma(TM) \}, \quad 
\Lambda^3_{27}
=\{\gamma\in\Lambda^3:
\gamma\wedge\varphi=0,\
\gamma\wedge *\varphi=0\}.
\]
The intrinsic torsion of the $\G_2$-structure, usually identified with $\nabla^{LC}\varphi$, can be encoded by four differential forms
$\tau_0,\tau_1,\tau_2,\tau_3$ determined by:
\begin{align}
 \dd\varphi&=\tau_0*\varphi+3\tau_1\wedge\varphi+*\tau_3,\label{eq:torsion1}\\
 \dd*\varphi&=4\tau_1\wedge*\varphi+\tau_2\wedge\varphi,
 \label{eq:torsion2}
\end{align}
where
$ \tau_0\in C^\infty(M),
 \tau_1\in\Omega^1(M),
 \tau_2\in\Omega^2_{14}(M),
 \tau_3\in\Omega^3_{27}(M).
$
This above representation-theoretic description is due to Bryant \cites{Bryant2006} and is naturally related to the Fern\'andez--Gray classification of the intrinsic torsion \cite{FernandezGray1982}, see also \cite{ChiossiSalamon2002}.

Special classes of $\G_2$-structures can be obtained by imposing simple conditions on the torsion forms. Some of the most important classes are:
\begin{itemize}
    \item \emph{nearly parallel} $\G_2$-structures, characterised by
    $
    \dd\varphi=\tau_0 *\varphi$,  where $\tau_0\neq 0$;
    
    \item \emph{closed} $\G_2$-structures, characterised by $
    \dd\varphi=0,
    $ or
    equivalently, $\tau_0=\tau_1=\tau_3=0$;
    \item \emph{co-closed} $\G_2$-structures, characterised by
    $
    \dd *\varphi=0,
    $ or
    equivalently, $\tau_1=\tau_2=0$;
    \item \emph{purely co-closed} $\G_2$-structures, which are co-closed and also satisfy $
    \tau_0=0
    $;
    \item \emph{$\G_2$-structures with torsion}, also called `\emph{integrable}', characterised by
    $
    \tau_2=0.
    $
\end{itemize}
An important topological criterion in the \textit{compact} torsion-free case is that the holonomy group is equal to $\G_2$ precisely when the fundamental group is finite. This is essentially a consequence of the Ricci flatness and Berger's classification \cite{Joyce2000}. By contrast there is considerable topological flexibility in the co-closed case: every compact spin $7$-manifold admits a co-closed $\G_2$-structure by the $h$-principle of Crowley--Nordstr\"om \cite{CrowleyNordstrom2015}. On the other hand, the topological existence problem for closed \(G_2\)-structures is much less understood. An important open problem is whether $S^7$ can admit a closed (and hence exact) $\G_2$-structure. Our main interest in this article lies in $\G_2$-structures with torsion, which turn out can be characterised by the existence of a $\G_2$ connection with skew-symmetric torsion, see Theorem \ref{FriedrichIvanovTheorem} below.

\section{Metric connections with skew torsion and \texorpdfstring{$\G_2T$}{}-structures}

Let $(M^n,g)$ be a Riemannian manifold and let $\nabla$ be a metric
connection. Writing
\[
 \nabla_XY=\nabla_X^{LC}Y+A(X,Y),
\]
the metric condition is equivalent to saying that $g(A(X,Y),Z)$ is skew-symmetric in $Y,Z$. Suppose now that the torsion of
$\nabla$ is a $3$-form $T$, then we can write
\begin{equation}\label{eq:metric-connection}
 g(\nabla_XY,Z)=g(\nabla_X^{LC}Y,Z)+\frac12T(X,Y,Z).
\end{equation}
Connections with skew-symmetric torsion preserve the unparametrised geodesics of $g$ and appear
naturally on compact Lie groups and naturally reductive homogeneous spaces,
see \cites{Agricola2006,AgricolaFriedrich2004} for a broad survey. In other words, they have the same geodesics as the Levi-Civita connection. The Ricci tensor of \eqref{eq:metric-connection} can be expressed as
\begin{equation}\label{eq:ricci-torsion}
 \Ric^\nabla=\Ric^g-\frac14T_g^2-\frac12\delta_gT, 
 \end{equation}
 where
$(T_g^2)_{ij}:=T_{ipq}T_{jpq}$ is a symmetric $2$-tensor.  Thus, the skew-symmetric part of $\Ric^\nabla$ is controlled by $\delta_gT$ while its symmetric part is determined by $\Ric^g$ and a quadratic term in $T$.

We recall the two important classes of examples worth keeping in mind. Let $G/H$ denote a naturally reductive homogeneous space with reductive decomposition
$\mathfrak g=\mathfrak h\oplus\mathfrak m$. Its canonical connection has
torsion
\[
 T(X,Y,Z)=-\langle [X,Y]_{\mathfrak m},Z\rangle.
\]
The latter is skew-symmetric and both the torsion and curvature tensors are parallel with respect to the canonical connection. For a compact Lie group with a bi-invariant metric the one-parameter family of connections
\[
 \nabla_X^\lambda Y=\lambda[X,Y]
\]
has torsion $T^{\nabla^{\lambda}}(X,Y)=(2\lambda-1)[X,Y]$ and curvature
\[
 R^{\nabla^\lambda}(X,Y)Z=\lambda(1-\lambda)[Z,[X,Y]].
\]
The two values $\lambda=0,1$ give the flat Cartan--Schouten connections. This set of 
examples are basic testing grounds for statements involving closed torsion, characteristic flatness and generalized Ricci flow.

For a general $\G_2$-structure there need not exist a $\G_2$-connection whose torsion is a $3$-form. The obstruction corresponds precisely to the $\Lambda^2_{14}$-component of the intrinsic torsion:
\begin{theorem}[Friedrich--Ivanov \cite{FriedrichIvanov2002}]\label{FriedrichIvanovTheorem}
A $\G_2$-structure $\varphi$ admits a metric connection $\nabla^T$ with totally
skew-symmetric torsion such that
\[
 \nabla^T\varphi=0
\]
if and only if $\tau_2=0$. When it exists, this connection is unique and its
torsion is given by
\begin{equation}\label{eq:char-torsion}
 T_\varphi=\frac16\tau_0\varphi+*(\tau_1\wedge\varphi)-\tau_3.
\end{equation}
\end{theorem}

We shall call such a structure a $\G_2T$-structure. It is also common to call
it an \emph{integrable} $\G_2$-structure, but this rather misleading because typically an \emph{integrable $H$-structure} means it is torsion free, so we shall try to avoid this terminology. From \eqref{eq:torsion2}, the
condition $\tau_2=0$ is equivalent to
\begin{equation}\label{eq:lee}
 \dd*\varphi=\theta\wedge*\varphi,
 \quad \theta=4\tau_1,
\end{equation}
where $\theta$ is the Lee form. The characteristic connection is the natural
seven-dimensional analogue of the Bismut connection in Hermitian geometry. The conditions $\tau_2=0$ and \eqref{eq:char-torsion} provide the link between intrinsic $\G_2$ geometry and the Killing spinor equations that occur in string theory, which is what we discuss next.

\section{The heterotic \texorpdfstring{$\G_2$}{} system}

\subsection{From Killing spinors to geometry}

In heterotic supergravity one studies ten-dimensional backgrounds which,
after compactification, carry a Riemannian metric, a dilaton, a gauge connection and a $3$-form flux $H$. Supersymmetry requires a non-zero spinor
$\eta$ satisfying a gravitino equation with respect to a metric connection with skew torsion, together with a dilatino equation and a gaugino equation.
In dimension seven, the stabiliser of a non-zero real spinor is $\G_2$, so the system naturally produces a $\G_2T$-structure.

At the differential form level, one convenient normalisation of the supersymmetry equations is
\begin{equation}\label{eq:susy-g2}
 \tau_2=0,\quad
 \tau_1=\frac12 d\mu,
 \quad
 F_A\wedge*\varphi= R_{\Theta}\wedge*\varphi=0,
\end{equation}
where $\mu$ is the dilaton, and $A$ is a connection on an auxiliary gauge bundle $P \to (M,\varphi)$ and $\Theta$ is a connection on $TM$.  The
instanton equation in \eqref{eq:susy-g2} says that the curvature 2-forms
$F_A$ and $ R_{\Theta}$ are of type $\Lambda^2_{14}$. 

The crucial extra equation is the anomaly cancellation, or Bianchi identity, given by:
\begin{equation}\label{eq:bianchi}
 d H=\frac{\alpha'}4\left(
 \operatorname{tr}(F_A\wedge F_A)-
 \operatorname{tr}(R_\Theta\wedge R_\Theta)\right).
\end{equation}
Supersymmetry identifies $H$ with the characteristic torsion $T_\varphi$. The numerical constants in
\eqref{eq:susy-g2}--\eqref{eq:bianchi} depend on the conventions for the
dilaton, trace and string parameter, but the geometric content is invariant:
$\varphi$ is integrable, the gauge and tangent-bundle connections are
$\G_2$-instantons, and the failure of $T_\varphi$ to be closed is prescribed by
a Pontryagin $4$-form. This formulation is rooted in the work of Strominger
and in the spinorial description of Friedrich--Ivanov
\cites{Strominger1986,FriedrichIvanov2002}.
\begin{remark}
In general, the supersymmetric compactification of heterotic supergravity on $M$ yields a three dimensional manifold, either Minkowski or anti-de Sitter spacetime. This depends on whether the cosmological constant $\lambda=\tfrac{7}{12}\tau_0$ vanishes. In fact one usually includes the condition $d\tau_0=0$ as part of  the heterotic $\G_2$-system, but it was recently shown that \eqref{eq:susy-g2} and \eqref{eq:bianchi} already imply that $\tau_0$ is constant \cite{CoupledG2Instantons}. 
\end{remark}
The simplest solutions are torsion-free $\G_2$ metrics. Compact solutions with
non-zero torsion and constant dilaton were constructed on nilmanifolds in
\cite{FIUV2009}. Further constructions exploit fibrations and gauge fields;
for example, Clarke--Garcia-Fernandez--Tipler obtained compact solutions on
torus bundles over K3 surfaces and studied the infinitesimal moduli
\cite{ClarkeGFT2022}, see also \cite{dBarcoMorenoFow26} and the references therein.

A useful limiting picture is obtained by formally sending the string parameter
$\alpha'$ to zero.  The right-hand side of \eqref{eq:bianchi} then disappears and
one is led to the closedness of the $3$-form $T_{\varphi}$. This is precisely the strong $\G_2T$ condition that will be discussed in more detail below. Thus, strong $\G_2T$ geometry can be viewed as an intrinsic analogue of SKT geometry and also the zero-anomaly limit of the heterotic system.
The analogy is especially close at the level of the distinguished
connections: the characteristic $\G_2$ connection plays the role of the Bismut connection, while the closed $3$-form torsion plays the role of the pluriclosed torsion $d^c\omega$.

\subsection{A differential complex}

An important feature of $\G_2T$ structures  is the existence of a
canonical differential complex. Let $\pi_7$ and $\pi_1$ denote the projections
onto $\Omega^2_7$ and $\Omega^3_1$, respectively. One considers
\[
 0\longrightarrow\Omega^0
 \overset{d}{\longrightarrow}\Omega^1
 \overset{\check d}{\longrightarrow}\Omega^2_7
 \overset{\check d}{\longrightarrow}\Omega^3_1
 \longrightarrow0,
\]
where the differentials are obtained by composing the exterior differential $\dd$ with the appropriate
$\G_2$-projections. The identity $\check d^{\,2}=0$ is equivalent to
$\tau_2=0$ \cite{FernandezUgarte1998}. If $A$ is a connection on a vector
bundle, the analogous covariant complex satisfies
\[
 \check d_A^{\,2}=0
 \quad\Longleftrightarrow\quad
 \tau_2=0\ \text{ and }\ F_A\wedge*\varphi=0.
\]
This complex plays an important role in the study of the infinitesimal deformation of the heterotic system. De la Ossa--Larfors--Svanes introduced a differential operator coupling the
geometry of $TM$ to the gauge and tangent-bundle instantons and showed that the square of the projected operator vanishes precisely when the heterotic
$\G_2$ equations, including the Bianchi identity, hold
\cite{deLaOssaLarforsSvanes2018}. The corresponding first cohomology group controls infinitesimal deformations. This gives a seven-dimensional analogue of deformation complexes familiar from complex and gauge geometry.

\subsection{Generalized geometry and coupled \texorpdfstring{$\G_2$}{}-instantons}

The Bianchi identity \eqref{eq:bianchi} has a natural interpretation in generalized geometry.
To keep the discussion short, suppose that the gauge and tangent bundles data
have been combined into a single principal bundle $P$ with invariant pairing on its
Lie algebra. One then considers the transitive Courant, or string, algebroid
\[
 E=TM\oplus\ad P\oplus T^*M.
\]
Its Dorfman bracket is twisted by $H$ and by the curvature of the combined connection $\theta$
on $P$. The condition
\begin{equation}\label{eq:string-bianchi}
 d H=\langle F_\theta\wedge F_\theta\rangle_P
\end{equation}
is exactly what is needed for the bracket to satisfy the Courant identities. In this way, the Bianchi identity can be viewed as the Jacobi-type identity which makes the transitive Courant algebroid well-defined.

A Riemannian metric on $M$, together with the connection and a choice of
splitting, determines a generalized metric on $E$. Generalized Ricci tensors are defined after choosing, in addition, a
divergence operator. The latter records the dilaton/conformal information.
For heterotic $\G_2$ systems the resulting generalized Ricci tensor vanishes;
this relates the supersymmetry equations to generalized Ricci-flatness
\cites{GarciaFernandez2019,CoupledG2Instantons}. 

There is also a useful gauge-theoretic reformulation. If
$\varphi$ is a $\G_2T$-structure, let
\[
 \nabla_H^{\pm}=\nabla^{LC}\pm\frac12g_\varphi^{-1}H.
\]
For a connection $\theta$ on $P$, its curvature gives off-diagonal bundle maps
$\mathbf F$ and $\mathbf F^\dagger$, and one can form a connection on
$TM\oplus\ad P$ of the schematic form
\begin{equation}\label{eq:coupled-connection}
 \mathbf D_{(\varphi,H,\theta)}=
 \begin{pmatrix}
 \nabla_H^- & \mathbf F^\dagger\\
 -\mathbf F & d_\theta
 \end{pmatrix}.
\end{equation}
A triple $(\varphi,H,\theta)$ satisfying \eqref{eq:string-bianchi} is called
a \emph{coupled $\G_2$-instanton} when the curvature of
\eqref{eq:coupled-connection} is a $\G_2$-instanton. In particular, heterotic $\G_2$ systems provide such coupled instantons. Thus, this formulation allows one to combine the $\G_2$-instanton, generalized Ricci-flatness and Killing spinor equations into a single geometric structure, see \cite{CoupledG2Instantons} for further detail.
\section{Strong \texorpdfstring{$\G_2$}{}-structures with torsion}

We now return to the strong condition and discuss its geometric consequences.
\begin{definition}
A $\G_2T$-structure $\varphi$ is \emph{strong} if its characteristic torsion is
closed:
\begin{equation}\label{eq:strong}
 d T_\varphi=0.
\end{equation}
\end{definition}

When the right-hand side of the heterotic Bianchi identity \eqref{eq:bianchi} vanishes, this  reduces precisely to \eqref{eq:strong}. Thus, strong $\G_2T$ geometry can
be viewed both as a natural intrinsic problem and as the zero-anomaly limit of
the heterotic equations. It is also the $\G_2$ counterpart of the SKT condition
$dd^c\omega=0$ in Hermitian geometry.

This condition is substantially stronger than $\tau_2=0$. Indeed,
substituting \eqref{eq:char-torsion} into \eqref{eq:strong} yields a nonlinear
system coupling $\tau_0,\tau_1$ and $\tau_3$. Nevertheless, strong structures
need not be parallel, homogeneous, or characteristic Ricci-flat, see below for explicit examples.

The closure of the characteristic torsion has immediate consequences for the
intrinsic torsion forms.  In particular, the scalar component of $dT_\varphi$
provides a useful scalar-curvature identity.  With the normalisation
\eqref{eq:torsion1}-\eqref{eq:torsion2}, a strong $\G_2T$-structure satisfies
\begin{equation}\label{eq:scalar-strong}
 \operatorname{Scal}(g_\varphi)
 =10\,\delta\tau_1+\frac{49}{24}\tau_0^2+24|\tau_1|^2.
\end{equation}
After integration on a compact manifold, the divergence term drops out.  Thus,
the scalar curvature is constrained by the scalar torsion and the Lee form,
and this gives useful obstructions in the invariant and homogeneous settings \cite{FinoMartinMerchanRaffero2024}.

When, in addition, the characteristic Ricci tensor vanishes, formula
\eqref{eq:ricci-torsion} implies
\[
 \delta T_\varphi=0,\quad \operatorname{Ric}^{g_\varphi}=\frac14(T_\varphi)^2.
\]
Hence the characteristic torsion is harmonic.  Notice that the converse is false in general; indeed one of the important results of \cite{FinoFowdar2026} is precisely the construction of strong $\G_2T$-structures with closed and co-closed torsion but non-vanishing characteristic Ricci tensor.

\subsection{The twisted \texorpdfstring{$\G_2$}{} equations and homogeneous rigidity}

A useful formulation, motivated by generalized geometry, combines
\eqref{eq:lee} and \eqref{eq:strong}. In a common normalisation one considers
\begin{equation}\label{eq:twisted-g2}
 d\varphi\wedge\varphi=0,\quad
 d*\varphi=\theta\wedge*\varphi,\quad
 d\theta=0,\quad
 d T_\varphi=0.
\end{equation}
These are often referred to as the twisted $\G_2$ equations. They arise as a
seven-dimensional analogue of pluriclosed geometry and are naturally related
to generalized Ricci flow. However, despite the apparent flexibility of \eqref{eq:twisted-g2}, the homogeneous
compact case is very rigid.

\begin{theorem}[
\cite{FinoMartinMerchanRaffero2024}]
Let $M^7$ be a compact homogeneous manifold carrying an invariant strong
$\G_2T$-structure (but not torsion free). Then, up to finite cover, $M$ is diffeomorphic to
\[
 S^3\times T^4\quad\text{or}\quad S^3\times S^3\times S^1.
\]
Moreover, both types carry invariant solutions of the twisted $\G_2$ equations,
and compact non-flat solvmanifolds admit no invariant strong $\G_2T$-structure.
\end{theorem}
The above result suggests that to find genuinely new strong $\G_2T$ geometries, one should consider less symmetric situations.

\subsection{Characteristic Ricci-flatness}

For a strong structure, formula \eqref{eq:ricci-torsion} shows that
characteristic Ricci-flatness is equivalent to a very rigid interaction
between the torsion and the Lee form. A recent result gives a particularly
clean characterisation.

\begin{theorem}[\cite{FinoFowdar2026}]\label{thm: finofow}
Let $(M^7,\varphi)$ be a strong $\G_2T$-manifold and let $\nabla^T$ be its
characteristic connection. The following conditions are equivalent:
\begin{enumerate}[label=(\roman*)]
\item $\Ric^{\nabla^T}=0$;
\item $T_\varphi$ is co-closed and the vector field $\theta^\sharp$ preserves
      the $\G_2$-structure, that is
      $\mathcal L_{\theta^\sharp}\varphi=0$;
\item the Lee form is parallel with respect to the characteristic connection,
      $\nabla^T\theta=0$.
\end{enumerate}
\end{theorem}

No compactness assumption is required for the equivalence. On a compact
manifold it has immediate topological consequences; in particular, a compact
strong $\G_2T$-manifold with $\Ric^{\nabla^T}=0$ has non-trivial third
cohomology. This should be contrasted with the general strong case. There exist local
strong $\G_2T$ structures whose characteristic Ricci tensor does not vanish
\cite{FinoFowdar2026}. Thus, closed torsion by itself does not force the
characteristic connection to be Ricci-flat. These examples separate two
conditions that are frequently identified in the most symmetric constructions.
\begin{remark}
The condition of parallel characteristic torsion
\[
\nabla^T T_\varphi=0
\]
is very restrictive, independently of the strong condition \(dT_\varphi=0\). In  \cite{CleytonMoroianuSemmelmann2021} Cleyton--Moroianu--Semmelmann studied the local structure of Riemannian manifolds carrying a metric connection with parallel skew-symmetric torsion. They showed that such manifolds admit a canonical splitting of the tangent bundle and an associated Riemannian submersion, whose fibres are locally naturally reductive. Their results give a classification in several cases and reduce the general problem to a rather constrained class of geometries. In particular, the irreducible geometries which occur include naturally reductive spaces, nearly K\"ahler $6$-manifolds and nearly parallel $\G_2$-manifolds.

For a $\G_2T$-structure, one may therefore regard
\[
dT_\varphi=0,\qquad
\delta T_\varphi=0,\qquad
\nabla^T T_\varphi=0
\]
as three different conditions on the characteristic torsion. They are not equivalent in general, with parallel torsion being the most restrictive from the local geometric point of view. In the strong $\G_2T$ setting considered in \cite{FinoFowdar2026}, additional genericity assumptions on the intrinsic torsion imply
\[
\nabla^T R^T=0,
\]
which implies naturally reductive. In the circle reduction considered below, a further genericity assumption forces \(R^T=0\).
\end{remark}
\subsection{The  relation with \texorpdfstring{$\SU(3)$}{}-structures}

A particularly effective way of constructing and studying \(\G_2\)-structures is to start with \(\SU(3)\)-structures in dimension six.  Let $(N^6,\omega,\Omega)$ be endowed with
an $\SU(3)$-structure, where
\[
 \Omega=\Omega^++i\Omega^-,\quad
 \omega\wedge\Omega^+=0,\quad
 \frac{2}{3}\omega^3=\Omega^+\wedge\Omega^-.
\]
On an interval $I$ with coordinate $t$, the $3$-form defined by
\begin{equation}\label{eq:su3-to-g2}
 \varphi=\omega\wedge dt+\Omega^+
\end{equation}
determines a $\G_2$-structure on $I\times N$, with
\[
 *\varphi=\frac12\omega^2-dt\wedge\Omega^-.
\]
Conversely, an oriented hypersurface orthogonal to a unit vector in a manifold endowed with a
$\G_2$-structure inherits an $\SU(3)$-structure.  This was developed by Hitchin \cite{Hitchin2001} as a way of constructing torsion free $\G_2$-structures starting from \emph{half flat $\SU(3)$-structures}, see also \cite{Bryant2010}. This elementary observation is particularly powerful when a symmetry reduces a seven-dimensional problem to an ODE on a family of homogeneous six-manifolds, see for instance \cite{MadsenSalamon2013}. In fact the original Bryant-Salamon examples can be obtained in this way.

The intrinsic torsion of an $\SU(3)$-structure is uniquely determined by two functions $\pi_0,\sigma_0$, two $1$-forms $\nu_1,\pi_1$, two primitive real $(1,1)$-forms $\pi_2,\sigma_2$, and a primitive real $3$-form $\nu_3$ of type $(2,1)+(1,2)$, see \cites{ChiossiSalamon2002,BedulliVezzoni2007}. With the conventions adopted here, these torsion forms are defined by the identities
\begin{align*}
 d\omega&=-\frac32\sigma_0\Omega^+
           +\frac32\pi_0\Omega^-+\nu_1\wedge\omega+\nu_3,\\
 d\Omega^+&=\pi_0\omega^2+\pi_1\wedge\Omega^+-\pi_2\wedge\omega,\\
 d\Omega^-&=\sigma_0\omega^2+J\pi_1\wedge\Omega^+-\sigma_2\wedge\omega.
\end{align*}
Thus, half-flat structures are characterised by
$d\Omega^+=0$ and $\omega \w d\omega=0$.  Next we discuss the use of $\SU(3)$-structures in the study of strong $\G_2T$-structures.

\subsection{Reduction to six dimensions}

We now assume that $\theta$ is non-zero and $\nabla^T$-parallel. In view of Theorem \ref{thm: finofow} above, it follows that $\theta^\sharp$ is Killing and preserves the $\G_2$-structure. It is therefore natural to consider the geometry transverse to the flow of \(\theta^\sharp\). Locally, and globally under the usual regularity assumptions on this flow, the resulting six-dimensional quotient carries an \(\SU(3)\)-structure. Moreover, the strong condition and the characteristic Ricci-flat equations can be expressed entirely in terms of the intrinsic torsion of this \(\SU(3)\)-structure and the curvature of the corresponding circle connection.

This construction was developed in \cite{FinoFowdar2026}*{\S 6}, which we briefly recall here. More precisely, after normalising the Lee form, one may write
\begin{equation}\label{eq:reduction-su3}
 \varphi=\eta\wedge\omega+\Omega^+,
 \quad
 *\varphi=\frac12\omega^2-\eta\wedge\Omega^-,
\end{equation}
where \(\eta\) is the connection \(1\)-form associated with the Lee direction and its curvature $F=d\eta$ is a primitive \((1,1)\)-form on the quotient. After a constant phase rotation of $\Omega$, denoted by $\Psi$, the induced
$\SU(3)$-structure is half-flat and satisfies
\begin{equation}\label{eq:cocoupled-su3}
 d\omega\wedge\omega=0,\quad
 d\Psi^+=0,\quad
 d\Psi^-=\kappa\,\omega^2,
\end{equation}
for a constant $\kappa>0$ determined by $\tau_0$. In particular the quotient is co-coupled; it is not complex in the non-trivial case, and its Nijenhuis
tensor is a non-zero $3$-form.

The converse provides a Gibbons--Hawking type construction of such \(\G_2T\)-structures from six-dimensional data. Suppose that
$(P^6,\omega,\Psi^\pm)$ satisfies \eqref{eq:cocoupled-su3}, and let
\(F\in\Lambda^{1,1}_0(P)\) be a closed integral \(2\)-form. The strong condition is then equivalent to the Bianchi-type equation
\begin{equation}\label{eq: reduced bianchi}
\Delta_\omega\omega
=4\kappa^2\omega+*_\omega(F\wedge F).
\end{equation}
The circle bundle with curvature \(F\) then carries a characteristic Ricci-flat strong \(\G_2T\)-structure. Conversely, every such structure with non-zero Lee form is locally obtained in this way \cite{FinoFowdar2026}.

Thus, the condition of parallel non-zero Lee form allows us to completely encode the problem of finding characteristic Ricci flat strong $\G_2T$-structure as a problem on a 6-manifold with suitable $\SU(3)$-structure. It is worth pointing out that the usual $\SU(3)$ heterotic system is formulated on a Hermitian threefold and is commonly referred to as the Hull--Strominger system. In contrast, the quotient arising here is not Hermitian: its almost complex structure is non-integrable, although its Nijenhuis tensor is totally skew-symmetric.

\section{Relation to generalized Ricci flow}

Let $g$ be a Riemannian metric and $H$ a closed $3$-form. After a choice of gauge, the generalized Ricci flow is given by 
\begin{equation}\label{eq:grf}
\frac{\partial g}{\partial t}
=-2\Ric^g+\frac12H_g^2,
\qquad
\frac{\partial H}{\partial t}=-\dd\delta_gH,
\end{equation}
where $H$ is assumed to be initially closed and stays in its cohomology class. The stationary points correspond to generalized Ricci flat metrics with closed torsion, see \eqref{eq:ricci-torsion}. In the Hermitian case it is known that \eqref{eq:grf} is equivalent to Streets-Tian pluriclosed flow \cite{StreetsTian2010} up to diffeomorphism. This naturally motivates the question whether such an analogue also exists for  $\G_2$-structures; this would give a parabolic method of finding characteristic Ricci flat strong $\G_2T$-structures. We investigated precisely this problem in \cite{FinoFowdar2026}. 

Let
$
\psi=*_\varphi\varphi
$
and denote by $\diamond$ the infinitesimal $\mathrm{GL}(7)$-action on differential forms. Starting from a strong $\G_2T$-structure, consider
\begin{equation}\label{eq:g2-generalized-ricci-flow}
\frac{\partial\psi}{\partial t}
=
\left(-\Ric(g_\varphi)+\frac14T_\varphi^2\right)\diamond\psi
-4\lambda\mathcal L_{\tau_1^\sharp}\psi
+S_7^2\diamond\psi,
\end{equation}
where $\lambda\in\mathbb R$ and $S_7^2\in\Omega^2_7(M)$.
The term $S_7^2\diamond\psi$ lies in the $\Lambda^4_7$-component of the variation and as such it only varies the $\G_2$-structure but  not the underlying metric.

Using the general variation formulas for a $\G_2$-structure, together with $\theta=4\tau_1$, one finds that \eqref{eq:g2-generalized-ricci-flow} induces
\begin{equation}\label{eq:induced-generalized-ricci-flow}
\frac{\partial g_\varphi}{\partial t}
=
-2\Ric(g_\varphi)
+\frac12T_\varphi^2
-\lambda\mathcal L_{\theta^\sharp}g_\varphi.
\end{equation}
Thus, when $H=T_\varphi$, the induced metric satisfies the first equation of \eqref{eq:grf} up to the diffeomorphism term generated by the Lee vector field. In particular, for $\lambda=0$ it agrees exactly with the metric equation in generalized Ricci flow. Conversely, every $\G_2$-flow inducing \eqref{eq:induced-generalized-ricci-flow} is necessarily of the form \eqref{eq:g2-generalized-ricci-flow}; the $\Lambda^2_7$-term records the part of the evolution of the $\G_2$-structure which is not seen by the metric. It is shown in \cite{FinoFowdar2026} that for suitable choices of $S_7^2$, the flow admits a unique smooth solution for short time on compact manifolds.

So the remaining question is whether \eqref{eq:g2-generalized-ricci-flow} can be chosen so that the torsion evolves by the second equation in \eqref{eq:grf}. At present, it is not known whether one can choose $S_7^2$ so that both
\[
\tau_2=0,\qquad \dd T_\varphi=0
\]
are preserved by the flow. More precisely, one would like to choose $S_7^2$ so that
\[
\frac{\partial T_\varphi}{\partial t}
=
-\Delta_\varphi T_\varphi
-4\lambda\mathcal L_{\tau_1^\sharp}T_\varphi,
\]
which is the gauge-fixed form of the evolution of $H$ in generalized Ricci flow. In that case, both the metric and the characteristic torsion would evolve by generalized Ricci flow, while the $\Lambda^2_7$-term would determine the remaining evolution of the $\G_2$-structure.

There is a second class of flows which is more naturally adapted to co-closed $\G_2$-structures. If $\psi$
is closed, one may consider
\begin{equation}\label{eq:modified-coclosed-flow}
\frac{\partial\psi}{\partial t}
=
-\dd T_\varphi
+2\dd\bigl((C-\lambda\tau_0)\varphi\bigr),
\end{equation}
where $C,\lambda\in\mathbb R$. For $C=\lambda=0$, the critical points are precisely the co-closed strong $\G_2T$-structures. Other choices of the parameters recover the Laplacian coflow and the modified Laplacian coflow. For an explicit open interval of values of $\lambda$, containing the modified Laplacian coflow, equation \eqref{eq:modified-coclosed-flow} admits a unique short-time solution after gauge fixing on compact manifolds \cite{FinoFowdar2026}.

The two flows therefore address different aspects of strong $\G_2T$ geometry. The first is designed so that the induced metric follows generalized Ricci flow, while the second preserves the co-closed condition and is closely related to the Laplacian coflow. A third approach is to instead assume that $\psi$ is only conformally closed. This allows more flexibility and leads to an extra choice of a conformal function. This approach was developed independently by
Garcia-Fernandez, Moreno, Payne and Streets \cite{GFMPS2025},
and by  Karigiannis, Picard and Suan \cite{KPS2026}. In the former, the fixed points are large-volume solutions of the heterotic $\G_2$ system with $\tau_0=0$, while the latter is motivated by a seven-dimensional lift of the anomaly flow on complex threefolds.

\addtocontents{toc}{\protect\setcounter{tocdepth}{2}}
\section{Examples}\label{sec: Examples}
\subsection{Review of known examples}\label{sec: old examples}

We begin this section by recalling the previously known examples of strong $\G_2T$-structures. A common feature of all these complete examples is that the underlying metric is always a Riemannian product and as such the holonomy group of the Levi-Civita is \textit{strictly} contained in $\SO(7)$.

\begin{example}
    [Homogeneous examples]
It was shown in \cite{FinoMartinMerchanRaffero2024} that there are exactly two compact homogeneous 7-manifolds, up to covering, that can carry a left-invariant strong $\G_2T$-structure: $S^3\times \mathbb{T}^4$ and $S^3\times S^3\times S^1$. 

$\bullet$ An explicit example on $S^3 \times \mathbb{T}^4$ can be described as follows. Consider the global co-framing $\{e^i\}_{i=1}^7$ satisfying the structure equations:
\begin{equation*}
de^1 = e^{23}, \quad 
de^2 = e^{31}, \quad 
de^3 = e^{12},  
\end{equation*}
and $de^i= 0$ for $i=4,5,6,7$. Then a strong $\G_2T$-structure is given by
\[
\vp = e^{145}+e^{167} +  e^{246}-e^{257}  +  e^{347}+e^{356}  - e^{123}.
\]
The induced metric $g_{\vp}$ is the standard product metric and has holonomy algebra $\mathfrak{hol}(\nabla^{LC})=\mathfrak{so}(3)\oplus\{0\}\subset\mathfrak{so}(3)\oplus\mathfrak{so}(4)\subset \mathfrak{so}(7).$

$\bullet$ An explicit example on $S^3 \times S^3 \times S^1$ can be described as follows. Consider the global co-framing $\{e^i\}_{i=1}^7$ satisfying the structure equations:
$de^1 = 0$ and  
\begin{equation*}
de^2 = 2\sqrt{2} e^{45}, \quad 
de^3 = 2\sqrt{2} e^{67}, \quad 
de^4 = 2\sqrt{2} e^{52}, \quad 
de^5 = 2\sqrt{2} e^{24}, \quad 
de^6 = 2\sqrt{2} e^{73}, \quad
de^7 = 2\sqrt{2} e^{36}.  
\end{equation*}
Then a strong $\G_2T$-structure is given by
\[
\vp =e^{123}+e^{145}+e^{167}+e^{246}-e^{257}-e^{347}-e^{356}.
\]
Since the induced metric is again a product metric, the holonomy algebra is  $\mathfrak{hol}(\nabla^{LC})=\mathfrak{so}(3)\oplus\mathfrak{so}(3)\subset\mathfrak{so}(6)\subset \mathfrak{so}(7).$

\end{example}

\begin{remark}
Note that the result in \cite{FinoMartinMerchanRaffero2024} does not assert that these are unique strong $\G_2T$-structures on these homogeneous spaces. Indeed there is another homogeneous strong $\G_2T$-structure on $S^3 \times S^3 \times S^1$ which is not diffeomorphic to the above example but nonetheless induces the same Riemannian metric. This second example will be described in Proposition \ref{prop:static} below. 
\end{remark}

\begin{example}
    [Inhomogeneous product examples]
The above example on $S^3 \times \mathbb{T}^4$ can be generalised to $S^3 \times M^4$, where $M$ is a hyperK\"ahler $4$-manifold. Denoting the (self-dual) hyperK\"ahler $2$-forms by $\omega_1,\omega_2,\omega_3$, we can define a strong $\G_2T$-structure by
\[
\vp = e^1 \w \omega_1 +  e^2 \w \omega_2  +  e^3 \w \omega_3 - e^{123}.
\]
Assuming $M$ has full holonomy $\SU(2)=\mathrm{Sp}(1)$, in this case we get $\mathfrak{hol}(\nabla^{LC})=\mathfrak{so}(3)\oplus\mathfrak{su}(2)\subset\mathfrak{so}(3)\oplus\mathfrak{so}(4)\subset \mathfrak{so}(7).$
\end{example}

\begin{example}
    [Local non-product examples]
In \cite{FinoFowdar2026} it is shown that there exist  cohomogeneity one strong $\G_2T$-structures with \textit{incomplete} metrics. One such explicit example is given as follows. Let $N^6$ denote a $2$-step nilmanifold with associated left-invariant co-framing satisfying
\[
    de^{2} = 0, \quad 
    de^{4} = 0, \quad
    de^{6} = 0, \quad
    de^{3} = e^{64}, \quad
    de^{5} = e^{26}, \quad
    de^{7} = e^{42}.
\]
We define on $\R^+_t\times N$ an orthonormal co-framing by
\begin{gather*}
 E^1=dt, \quad   E^2 = 2t^{1/2}e^2, \quad E^4 = 2t^{1/2}e^4, \quad E^6 = 2t^{1/2}e^6, \quad
    E^3 = e^3, \quad E^5 = e^{5}, \quad E^7 = e^{7},
\end{gather*}
and a $\G_2$ $3$-form by
\[
\vp=
E^{123}+E^{145}+E^{167}+E^{246}-E^{257}-E^{347}-E^{356}.
\]
It is shown in \cite{FinoFowdar2026} that this defines a strong $\G_2T$-structure. The associated metric is incomplete since it becomes singular at finite distance as $t\to 0$; nonetheless, it is complete for $t\to+\infty$. A direct calculation shows that $\mathfrak{hol}(\nabla^{LC})=\mathfrak{so}(7)$ and hence this metric is not locally a Riemannian product.
\end{example}

In \cite{FinoFowdar2026} we also showed that there exist several incomplete examples on open sets of $\R \times S^3 \times S^3$. These examples naturally raise the question: does there exists a \textit{complete} strong $\G_2T$-structure whose metric is not locally a product? We shall answer this question in the affirmative in the next section by constructing {complete} examples on $\R\times S^3\times S^3$.

\subsection{New cohomogeneity one examples}

We begin by recalling the set up in \cite{FinoFowdar2026}*{\S 7.3}. Consider the standard left invariant co-framing $\{e^i\}_{i=1}^6$ on $S^3 \times S^3$ satisfying
\[
de^1 = -2 e^{23}, \quad 
de^2 = -2 e^{31}, \quad 
de^3 = -2 e^{12}, \quad 
de^4 = -2 e^{56}, \quad 
de^5 = -2 e^{64}, \quad 
de^6 = -2 e^{45}. 
\]
By mixing the co-framing on the two $S^3$ factors,  we can define a new co-framing by
\begin{equation*}
\begin{alignedat}{3}
\eta^1 &{}= \frac{1}{2}(e^1-e^4), \qquad&
\eta^2 &{}= \frac{1}{2}(e^1+e^4),  \qquad&
\eta^3 &{}= \frac{1}{2}(e^2-e^5), \\
\eta^4 &{}= \frac{1}{2}(e^2+e^5), &
\eta^5 &{}= \frac{1}{2}(e^3-e^6), &
\eta^6 &{}= \frac{1}{2}(e^3+e^6).
\end{alignedat}
\end{equation*}
Using the above, we define an adapted $\G_2$ co-framing on $\R_t \times S^3 \times S^3$ by
\begin{equation}\label{eq:cohom1-coframe}
\begin{alignedat}{3}
E^1 &{}= 2B_1(t)\eta^1, \qquad&
E^3 &{}= 2B_2(t)\eta^3, \qquad&
E^5 &{}= 2B_2(t)\eta^5, \\
E^2 &{}= 2A_1(t)\eta^2, &
E^4 &{}= 2A_2(t)\eta^4, &
E^6 &{}= 2A_2(t)\eta^6.
\end{alignedat}
\end{equation}
and $E^7 = dt$, so that the $\SU(2)^2\times \mathrm{U}(1)$-invariant $\G_2$ $3$-form can be expressed as 
\begin{equation}\label{eq: vpansatz}
\vp:= E^{127}+E^{347}+E^{567} + E^{135}-E^{146}-E^{236}-E^{245}.
\end{equation}
The $\SU(2)^2$ symmetry is clear by the choice of left invariant $1$-forms, and the extra $\mathrm{U}(1)$ symmetry corresponds to the repeated appearance of the coefficient functions $A_2(t)$ and $B_2(t)$ in the above. This $\mathrm{U}(1)$ symmetry will be relevant in Section \ref{sec: quotient}.

Next we describe the torsion forms.
A direct calculation shows that the torsion form $\tau_0$ is always zero. To describe the remaining torsion forms, we first introduce the following functions:
\begin{align*}
 P &=\frac{A_1}{B_2^2}+\frac{2A_2}{B_1B_2}
 -2\frac{B_2'}{B_2}-\frac{B_1'}{B_1},\\
 Q &=\frac{A_1}{A_2^2}+\frac{B_1'}{B_1}
 +2\frac{A_2'}{A_2}-\frac{2B_2}{A_2B_1},\\
 R &=\frac{B_2'}{B_2}+\frac{A_2'}{A_2}
 -\frac{B_1}{A_2B_2}+\frac{A_1'}{A_1},\\
 M &=\frac{A_1'}{A_1}+\frac{A_2'}{A_2}
 +\frac{B_1'}{B_1}+\frac{B_2'}{B_2}
 -\frac{A_2}{B_1B_2}-\frac{B_2}{A_2B_1},\\
 N&= 2\frac{A_2'}{A_2}+2\frac{B_2'}{B_2}
 -\frac{2B_1}{A_2B_2}+\frac{A_1}{A_2^2}-\frac{A_1}{B_2^2}.
\end{align*}
The remaining $\G_2$ torsion forms can then be neatly expressed as follows:
\begin{gather}
    \tau_1 = \frac{N+2M}{12}E^7\\
    \tau_2 = \frac{N-M}{3}(-2E^{12}+E^{34}+E^{56})\\
    \tau_3 = (R-\frac{N}{4}-\frac{M}{2}) (E^{145}+E^{136}) + (Q-\frac{N}{4}-\frac{M}{2}) E^{235}+ (P+\frac{N}{4}+\frac{M}{2})E^{246}.
\end{gather}
Thus, we see that the $\G_2T$ condition is equivalent to setting $N=M$. Under this assumption, we have
\begin{equation}
    T_{\vp} = (N-R) (E^{145}+E^{136}) + (N-Q) E^{235}- (N+P)E^{246}.
\end{equation}
Since $T_{\vp}$ is of the form 
\begin{equation}
    T_{\vp} = f_1(t) (\eta^{145}+\eta^{136}) + f_2(t) \eta^{235}+ f_3(t)\eta^{246}.
\end{equation}
and from the structure equations we also have
\begin{equation*}
d(\eta^{145}+\eta^{136})=-4\eta^{3456},\quad
d\eta^{235}=2\eta^{3456}-2\eta^{1256}-2\eta^{1234},\quad
d\eta^{246}=2\eta^{3456}+2\eta^{1256}+2\eta^{1234},
\end{equation*}
it is easy to see that imposing the strong condition, i.e. $dT_{\vp}=0$, implies that 
\[
T_{\vp} = 4\kappa(e^{123}+e^{456}),
\]
where $\kappa$ is a constant. In other words, $ f_i'(t)=0$ and $f_1=f_2=f_3$; hence $f_1=f_2=f_3=8\kappa\in 
\R$.

Thus, the strong $\G_2T$ condition is equivalent to the system:
\begin{align}
 A_1'&=-\frac{(A_1^2-\kappa)(A_2^2-B_2^2)}{2A_2^2B_2^2},
 \label{eq:A1system}\\
 A_2'&=\frac12\left[
 \frac{\kappa+B_1^2+B_2^2-A_2^2}{B_1B_2}
 -\frac{\kappa+A_1^2}{A_1A_2}\right],\label{eq:A2system}\\
 B_1'&=\frac{\kappa+A_2^2-B_1^2+B_2^2}{A_2B_2},
 \label{eq:B1system}\\
 B_2'&=\frac12\left[
 \frac{\kappa+A_2^2+B_1^2-B_2^2}{A_2B_1}
 +\frac{\kappa+A_1^2}{A_1B_2}\right].\label{eq:B2system}
\end{align}
Note that when $\kappa=0$, we recover exactly the torsion free condition derived in \cite{Bogoyavlenskaya2013}*{Lemma 1}. In the latter it also shown that there exists a 1-parameter family of complete $\G_2$ metrics solving this system; this is often called the $\mathbb{B}_7$ family in the physics literature. In this case, the metrics can be extended smoothly to $\R^4 \times S^3$, which corresponds to the compactification of $\R\times S^3 \times S^3$ at one end by $S^3$. There are only two $\G_2$-structures in this $1$-parameter family which are known explicitly; the Bryant-Salamon one \cite{BryantSalamon1989} corresponding to: 
\[A_1=A_2= \frac{r}{3} \sqrt{1-r^{-3}}\sim \frac{r}{3}\ \text{  and  }\ B_1=B_2 = \frac{r}{\sqrt{3}},\]
where $E^7= \frac{dr}{\sqrt{1-r^{-3}}}$ and $r\in[1,+\infty)$, 
and the Brandhuber-Gomis-Gubser-Gukov (BGGG) one \cite{BrandhuberGomisGubserGukov2001} corresponding to: 
\begin{align*}
A_1=\frac{\sqrt{(r-9/4)(r+9/4)}}{\sqrt{(r-3/4)(r+3/4)}}\sim 1,\ \
A_2=\sqrt{\frac{(r-9/4)(r+3/4)}{3}}\sim \frac{r}{\sqrt{3}}\\
B_1=\frac{2r}{3},\ \
B_2=\sqrt{\frac{(r+9/4)(r-3/4)}{3}}\sim \frac{r}{\sqrt{3}}
\end{align*}
where $E^7=\frac{\sqrt{(r-3/4)(r+3/4)}}{\sqrt{(r-9/4)(r+9/4)}}dr$ and $r\in[\frac{9}{4},+\infty).$

The Bryant-Salamon metric is asymptotically conical (AC) i.e.
the metric is asymptotic the cone metric over a nearly K\"ahler $6$-manifold which in this case is $S^3 \times S^3$. On the other hand, the BGGG metric is asymptotically locally conical (ALC) i.e. it is asymptotic to an $S^1$ bundle of constant fibre length over a Calabi-Yau cone which is this case is over $S^2 \times S^3$. Aside from the Bryant-Salamon metric, all the torsion free metrics is the $1$-parameter family are ALC. We refer the reader to the introduction in \cite{FoscoloHaskinsNordstrom2021JEMS} which gives a comprehensive overview of all these cohomogeneity one torsion free $\G_2$-structures. 

Going back to the general system \eqref{eq:A1system}-\eqref{eq:B2system}, we first show that it has a stationary point:
\begin{proposition}\label{prop:static}
    There is, up to sign choices of the coefficient functions, a unique stationary solution to \eqref{eq:A1system}-\eqref{eq:B2system}, and it occurs precisely for $\kappa<0$. Moreover, in this case, $A_1^2=A_2^2=B_1^2=B_2^2=-\kappa$: this corresponds to the product metric on $\R \times S^3 \times S^3$ with the same bi-invariant metric on each $S^3$ factor.
\end{proposition}
\begin{proof}
Since we need a non-degenerate metric, we shall assume that $A_1A_2B_1B_2\neq0$. At a static solution all four derivatives vanish. From \(B_1'=0\) we obtain
\begin{equation}
 B_1^2=\kappa+A_2^2+B_2^2.
 \label{eq:B1relation}
\end{equation}
Substituting \eqref{eq:B1relation} into \(A_2'=0\) gives
\begin{equation}
 2A_1A_2(\kappa+B_2^2)
 =B_1B_2(\kappa+A_1^2).
 \label{eq:relation1}
\end{equation}
Similarly, substituting \eqref{eq:B1relation} into \(B_2'=0\) gives
\begin{equation}
 2A_1B_2(\kappa+A_2^2)
 =-B_1A_2(\kappa+A_1^2).
 \label{eq:relation2}
\end{equation}
Multiplying \eqref{eq:relation1} by \(A_2\) and
\eqref{eq:relation2} by \(B_2\), and adding the result eliminates the term with
\(\kappa+A_1^2\).  Since \(A_1\neq0\), this yields
\begin{equation}
 \kappa(A_2^2+B_2^2)+2A_2^2B_2^2=0.
 \label{eq:keyrelation}
\end{equation}
From the latter, it follows that $\kappa<0$ otherwise $A_2=B_2=0$.

On the other hand, \(A_1'=0\) gives
\begin{equation}
 (A_1^2-\kappa)(A_2^2-B_2^2)=0.
 \label{eq:A1static}
\end{equation}
Since \(\kappa<0\), we must have 
$ A_2^2=B_2^2.$ Equation \eqref{eq:keyrelation} then becomes
\[
 \kappa =-A_2^2.
\]
Substituting into \eqref{eq:B1relation} gives
\[
 B_1^2=-\kappa.
\]
Finally, \eqref{eq:relation1} becomes
\[
 0=B_1B_2(\kappa+A_1^2),
\]
and the result follows.
\end{proof}

\begin{remark}
Using the result in \cite{Fowdar2020thesis}*{Lemma 6.4.10}, one can check directly that any $\SU(2)^2 \times \mathrm{U}(1)$ strong $\G_2T$-structure given by \eqref{eq: vpansatz} which extends smoothly to the spinor bundle of $S^3$ is necessarily torsion free. In other words, the smooth extension condition forces $\kappa=0$ in \eqref{eq:A1system}-\eqref{eq:B2system}. This was already observed in \cite{FinoFowdar2026} in the case of $\SU(2)^3$-invariance.
\end{remark}
We shall next describe two complete non-torsion-free strong \(\G_2T\)-structures on \(\mathbb{R}\times S^3\times S^3\). In both cases, one end is asymptotic to the homogeneous stationary solution above, while the other end is asymptotic to a torsion-free \(\G_2\)-metric. The first has an AC end and the second an ALC end.

\subsubsection{\texorpdfstring{$\SU(2)^3$}{} example: AC end}
\begin{figure}[htp]
\centering
\includegraphics[width=.5\textwidth]{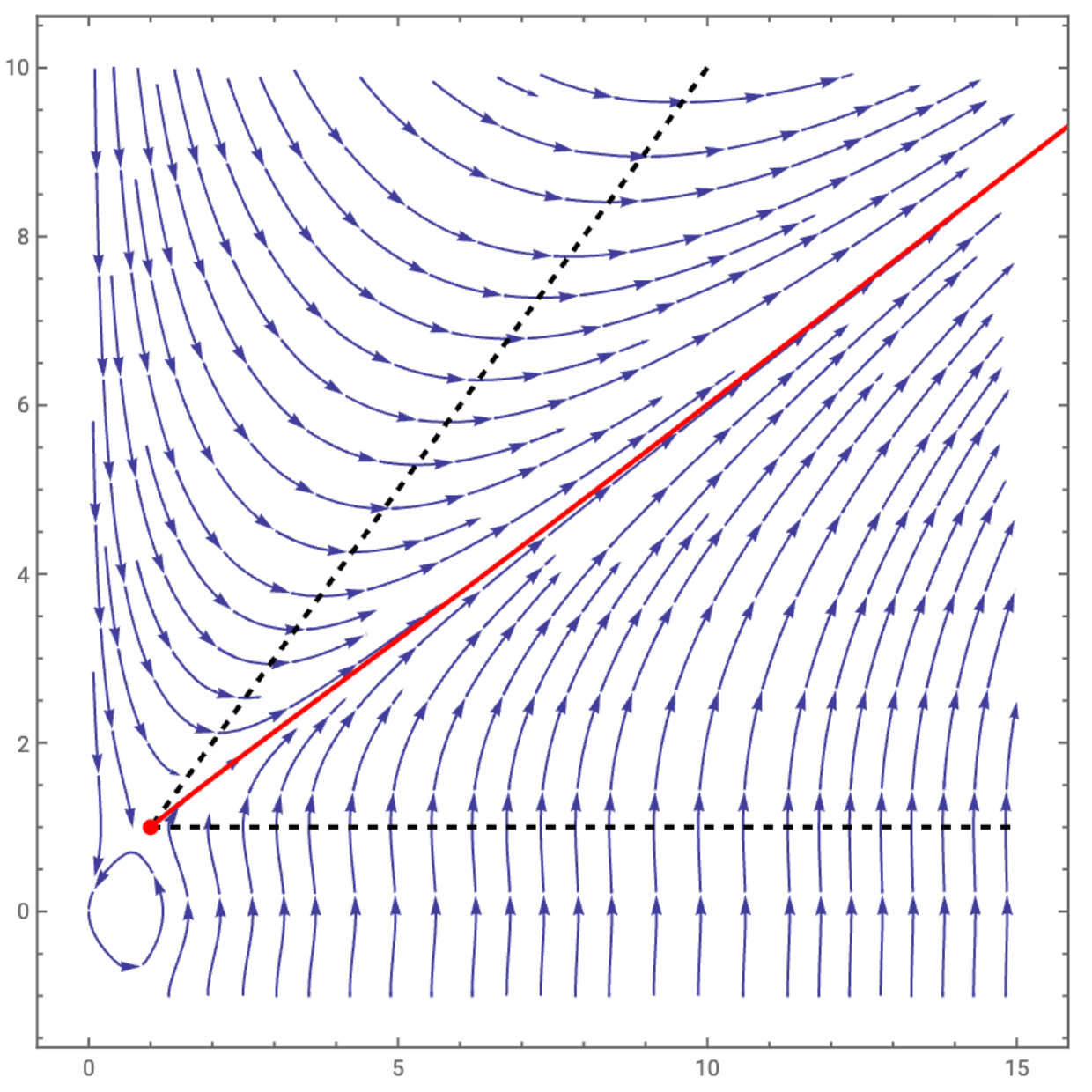}
\caption{{$(h,f)$-phase diagram of (\ref{equ: single ode}) for $\kappa=-\frac{1}{4}$.}}
\label{fig:figure1}
\end{figure}

We now specialise \eqref{eq:A1system}-\eqref{eq:B2system} to the case when $f=2A_1=2A_2$ and $h=2B_1=2B_2$. In other words, the $\G_2$-structure is now $\SU(2)^3$ invariant. This is precisely the situation considered in \cite{FinoFowdar2026}*{\S 7.3 (63)-(64)}:
\begin{equation}
 h'=\frac{8\kappa+2f^2}{hf},\qquad
 f'=\frac{(f^2-4\kappa)(h^2-f^2)}{h^2f^2}.\label{equ: single ode}
\end{equation}
In \cite{FinoFowdar2026}, we showed that none of the solutions to \eqref{equ: single ode} can be extended smoothly to $\R^4 \times S^3$ for $\kappa\in \R$,  except for $\kappa=0$ which corresponds to the torsion free case. However, it turns out there is one complete solution to \eqref{equ: single ode} which is instead well-defined on $\R\times S^3 \times S^3$; this occurs only for $\kappa<0$ as in Proposition \ref{prop:static}. Without loss of generality, in this case we can set $\kappa=-\frac{1}{4}$ so that the static solution occurs at $(h,f)=(1,1)$. The strong $\G_2T$ complete metric corresponds to a branch of the unstable manifold to this stationary point; this is depicted in red in Figure \ref{fig:figure1}. 
We shall now establish the existence of this solution rigorously and show that it has the same asymptotic behaviour as the Bryant-Salamon AC solution $A_1 \sim \frac{t}{3}$ and $B_2 \sim \frac{t}{\sqrt{3}}$ as $t\to +\infty$. 

\begin{proposition}\label{prop:su23branch} 
    The stationary point $(h,f)=(1,1)$ is a hyperbolic critical point with both the stable and unstable manifold $1$-dimensional. Moreover, one branch of the unstable manifold lies in the region
\[
 \mathcal D=\{(h,f):h>f>1\}
\]
and remains in $\mathcal D$ for all forward time.
\end{proposition}
\begin{proof}
    Computing the Jacobian of the system \eqref{equ: single ode}, with $\kappa=-\frac{1}{4}$, at $(h,f)=(1,1)$ gives:
    \[A=
\begin{pmatrix}0&4\\4&-4\end{pmatrix}.
\]
It is easy to check that the eigenvalues are given by $\lambda_\pm=-2\pm2\sqrt5\approx -6.47,2.47$. Since these are real with one positive and one negative value, the stable/unstable manifold theorem asserts that there are 1-dimensional stable and unstable manifolds at $(1,1)$. In particular, it follows that for the unstable manifold $(h(t),f(t))\to (1,1)$ as $t\to -\infty$; in fact the convergence is at an exponential rate. Observe also that there are 2 branches of the unstable manifold, one for which $(h(t),f(t))\to (0,0)$ as $t$ increases (one can show that this occurs in finite time and as such the resulting metric is singular, so we shall disregard this branch) and the second one for which $(h(t),f(t))$ is tangent  at $(1,1)$ to the line defined by $f-1=\frac{\lambda_+}{4}(h-1)$ since this corresponds to the eigenvector of the Jacobian i.e $$
A\begin{pmatrix}h-1\\f-1\end{pmatrix}=\lambda_+\begin{pmatrix}h-1\\f-1\end{pmatrix}.
$$ Thus, it is easy to see that the latter branch lies in the region $\mathcal{D}$ near $(1,1)$; this is also clear from Figure \ref{fig:figure1}. 

We claim that this branch stays in $\mathcal{D}$ for all $t$. On the boundary $f=1$ and $h>1$, from \eqref{equ: single ode}, we have $h'(t)=0$ and $f'(t)=\frac{2(h(t)^2-1)}{h(t)^2}>0$. So the curve cannot cross the line $f=1$ to $f<1$. On the other hand, on the boundary $h=f=C>1$ we have  $f'=0$ and $h'=2(C^2-1)/C^2>0$ i.e. $h'-f'>0$. So the curve cannot cross the line $f=h$ to $f>h$. Thus, the unstable branch (in red) coming out of $(1,1)$ stays inside the region $\mathcal{D}$.
\end{proof}

\begin{lemma}\label{lemma:ac}
The branch in Proposition~\ref{prop:su23branch} is defined for all
$t\in\R$, and $h(t),f(t)\to+\infty$ as $t\to+\infty$. Moreover, $h= \frac{2}{\sqrt{3}}t+O(1)$ and $f= \frac{2}{{3}}t+O(1)$ i.e. the $\G_2$-structure is asymptotic to the Bryant-Salamon cone metric.
\end{lemma}
\begin{proof}
By the stable/unstable manifold theorem, we already know that $(h(t),f(t))$ converges exponentially to $(1,1)$  as $t\to -\infty$. So it only remains to consider the case $t\to +\infty$. 
In the region $\mathcal D=\{(h,f):h>f>1\}$, from \eqref{equ: single ode} we have
\[
 0< h'=2\frac{f}{h}\left(1-\frac1{f^2}\right)<2,
 \qquad
 0< f'=\left(1+\frac1{f^2}\right)
 \left(1-\frac{f^2}{h^2}\right)<2.
\]
Thus, $h$ and $f$ are strictly increasing and they cannot blow up in finite forward time. We now show that they cannot be bounded above and hence $f,h\to +\infty$. Suppose $h$ is bounded above, then so is $f$, and say $(h(t),f(t))\to(H,F)$. Then from the system of ODEs we must have $h'(t),f'(t)\to 0$ which implies that $F=H=1$, but this is a contradiction since $f,h$ are increasing. Now suppose  $f$ is  bounded from above and $h\to +\infty$; say $f(t)\to F$. Then again we get $f'(t)\to 1+\frac{1}{F^2} > 0$ and hence a contradiction. Thus, we have established that both $h(t),f(t)\to +\infty$ and hence the induced metric is complete.

Now we turn to the growth rate. First we set
\[
q=\frac{f^2}{h^2} \qquad\text{and}\qquad x=\log(h).
\]
Using \eqref{equ: single ode}, a direct calculation shows 
\begin{equation}\label{eq:qeq}
 \frac{d q}{d x}
 =\frac{-3h^2q^2+(1+h^2)q+1}{f^2-1}.
\end{equation}
Since we are in the region $\mathcal{D}$, it is clear that the denominator ${f^2-1}>0$ cannot vanish. Suppose that $ \frac{d q}{d x}=0$ then
\[
-3q^2+(\frac{1}{h^2}+1)q+\frac{1}{h^2}=0.
\]
Solving for $q$ gives
\[
 q=\frac{h^2+1+\sqrt{(h^2+1)^2+12h^2}}{6h^2};
\]
Denoting this solution by $q_*(h)$, it is not hard to check that $q_*(h)$ is a strictly decreasing function for $h>1$ and $q_*(h)\to \frac{1}{3}$ as $h\to +\infty$. For $q(h)<q_*(h)$, we have $\frac{dq}{dx}>0$, and for $q(h)>q_*(h)$, we have $\frac{dq}{dx}<0$. We shall use the curve $q_*(h)$ as a barrier.

From the stable/unstable manifold theorem, we know that near $(1,1)$ the solution curve $(h(t),f(t))$ is of the form
\[
f(t)-1=\frac{\lambda_+}{4}(h(t)-1)+O((h(t)-1)^2).
\]
Thus, using the latter expression, we obtain the Taylor series expansion
\[
q(t)=1+2(\frac{\lambda_+}{4}-1)(h(t)-1)+O((h(t)-1)^2),
\]
where $2(\frac{\lambda_+}{4}-1)\approx-0.76$. For $q_*(h)$, at $(1,1)$ we instead have
\[
q_*(h)=1-(h-1)+O((h-1)^2).
\]
This shows that near $(1,1)$ we have $q(t)>q_*(h)$,
i.e. $\frac{dq}{dx}<0$ for the solution. We claim that this holds for all time. Suppose that the first point of contact is $x_0$ where we have $q(h(t))=q_*(h)$, 
then $\frac{dq}{dx}\big|_{x_0}=0$ and also $\frac{d}{dx}(q-q_*)\big|_{x_0}\leq 0$ since $q(x)-q_*(x)>0$ for $x<x_0$. But $q_*(h)$ is also decreasing as a function of $h$ and so we find $\frac{d}{dx}q_*=\frac{dh}{dx}\frac{d}{dh}q_*=e^x\frac{d}{dh}q_*<0$. This shows that $\frac{d}{dx}(q-q_*)\big|_{x_0}>0$ which is contradiction and as such the solution curve cannot cross the barrier $q_*(h)$. So $q(h(t))$ is decreasing and bounded below by $\frac{1}{3}$; suppose the limit of $q(t)$ is $L\geq \frac{1}{3}$. If $L>1/3$, then the right-hand side of \eqref{eq:qeq} tends to the negative
constant $1-3L$ as $x\to+\infty$, which contradicts the fact that $q(t)$ converges.
Therefore, $L=\frac{1}{3}$ i.e.
\begin{equation}\label{eq:su23ratio}
 \frac fh\longrightarrow\frac1{\sqrt3}.
\end{equation}
Substituting in \eqref{equ: single ode} and using $f,h\to+\infty$ then gives
\begin{equation}\label{eq:su23AC}
h'(t)\longrightarrow\frac2{\sqrt3},\qquad
 f'(t)\longrightarrow\frac23,\qquad
 h(t)\sim\frac{2t}{\sqrt3},\quad f(t)\sim\frac{2t}{3}.
\end{equation}
To show that the end is indeed asymptotically conical, we need to improve the above rate of decay. Writing $u=q-\frac{1}{3}$, equation \eqref{eq:qeq} can now be expressed as
\[
u_x+3u=
\frac{\frac43h^{-2}-2uh^{-2}}
{\frac13+u-h^{-2}}.
\]
Since we already know that $u\to 0$ and $h\to+\infty$, we have
\[
u_x+3u=O(e^{-2x}),
\]
where recall $h:=e^x$. Hence
\[
\frac{d}{dx}(e^{3x}u)=e^{3x}(u_x+3u)=O(e^{x}).
\]
Integrating the latter shows $e^{3x}u=O(e^{x})$ and hence
\[
u=O(e^{-2x})=O(h^{-2})=O(t^{-2}).
\]
Thus, we have the improved rate of decay
\[
q=\frac{1}{3}+O(t^{-2}).
\]
Substituting into \eqref{equ: single ode} as before now yields: 
\[
h'(t)=\frac2{\sqrt3}+O(t^{-2}),\qquad
f'(t)=\frac23+O(t^{-2}),
\qquad
h(t)=\frac{2t}{\sqrt3}+O(1),\qquad
f(t)=\frac{2t}{3}+O(1).
\]
Hence this implies that the metric is AC to the Bryant-Salamon cone metric:
\[
g_{\mathrm{AC}}=dt^2+\frac{4}{3}t^2\big((\eta^1)^2+(\eta^3)^2+(\eta^5)^2\big)+\frac{4}{9}t^2\big((\eta^2)^2+(\eta^4)^2+(\eta^6)^2\big)
\]
at the rate $\nu=-1$ since $|g_\varphi-g_{cone}|_{g_{cone}}=O(t^{-1})$. 
Moreover, using the fact that 
$\partial_{x}^k(u)=O(e^{-2x})$ for all $k\in \mathbb{N}_0$, it is not hard to see that one can repeat the above argument to obtain $|\nabla_{cone}^k(g_\varphi-g_{cone})|_{g_{cone}}=O(t^{-1-k})$.
\end{proof}

The above can be summarised as follows:
\begin{theorem}\label{thm:mainsu(2)3}
There exists a complete strong $\G_2T$-structure on
$\R\times S^3\times S^3$ whose associated metric is asymptotically
cylindrical at one end, asymptotic
to the product of the Euclidean metric on $\R$ with a bi-invariant
metric on $S^3\times S^3$, and is asymptotically conical at the other end to the
Bryant--Salamon cone metric. Moreover, its torsion $3$-form is independent of $t$ and is given explicitly by
\[
T_{\vp}=-(e^{123}+e^{456}),
\]
where $e^{123}$ and $e^{456}$ denote the usual volume forms on the two
$S^3$ factors. 
\end{theorem}

It is not hard to see that the Riemannian metric associated to this complete strong $\G_2T$-structure is not locally a product. Indeed using the structure equations \eqref{eq:cohom1-coframe} and \eqref{equ: single ode}, one can directly compute the curvature of the associated metric $g_{\varphi}$, say with the aid of Mathematica, and verify that $\dim(\mathfrak{hol}(\nabla^{LC}))=21$. Thus, unlike the previously known examples, see \S \ref{sec: old examples}, the cohomogeneity one solution given in Theorem \ref{thm:mainsu(2)3} is a genuinely new example with $\mathfrak{hol}(\nabla^{LC})=\mathfrak
{so}(7)$. Moreover, by direct calculation one can check the following hold:
\begin{enumerate}
    \item The torsion $3$-form $T_{\vp}$ is co-closed. In fact this holds for any strong $\G_2T$-structure satisfying \eqref{eq:A1system}-\eqref{eq:B2system}.
    \item The $\G_2$ Lee form $\theta=4\tau_1=-\frac{2(3f^2-h^2)}{h^2f^3}dt$ is not parallel with respect to $\nabla^T$. In view  \cite{FinoFowdar2026}*{Theorem 3.14}, this is equivalent to the fact that $\theta^\sharp$ does not preserve $\varphi$, or also equivalently, the Ricci tensor of $\nabla^T$ is non zero. One can also easily see  that the Lee vector field $\theta^\sharp$
    is not Killing.
\end{enumerate}

\subsubsection{\texorpdfstring{$\SU(2)^2\times \mathrm{U}(1)$}{} example: ALC end}

In this section we investigate a different ansatz to the previous one. Going back to the system \eqref{eq:A1system}-\eqref{eq:B2system}, we now set $A_2=B_2$. This choice is particularly natural since from \eqref{eq:A1system} it immediately forces $A_1=a$ to be constant, and  without loss of generality we shall set $a=1$. Moreover, from \eqref{eq:A2system} and \eqref{eq:B2system}, we get
\[
0=A_2'-B_2'=-\frac{1+\kappa}{A_1 A_2}
\]
and hence this shows that we need $\kappa=-1$ for the equations to be consistent. Under these assumptions, the system \eqref{eq:A1system}-\eqref{eq:B2system} reduces to the pair:
\begin{equation}\label{eq:hfsystemALC}
 h'=\frac{2f^2-h^2-1}{f^2}, \qquad  f'=\frac{h^2-1}{2hf},
\end{equation}
where $h=B_1$ and $f=A_2=B_2$. The normalisation choice $a=1$ was made to ensure that \eqref{eq:hfsystemALC} has a stationary point at
$(h,f)=(1,1)$ as before, described in Proposition~\ref{prop:static}.
\begin{figure}[htp]
\centering
\includegraphics[width=.5\textwidth]{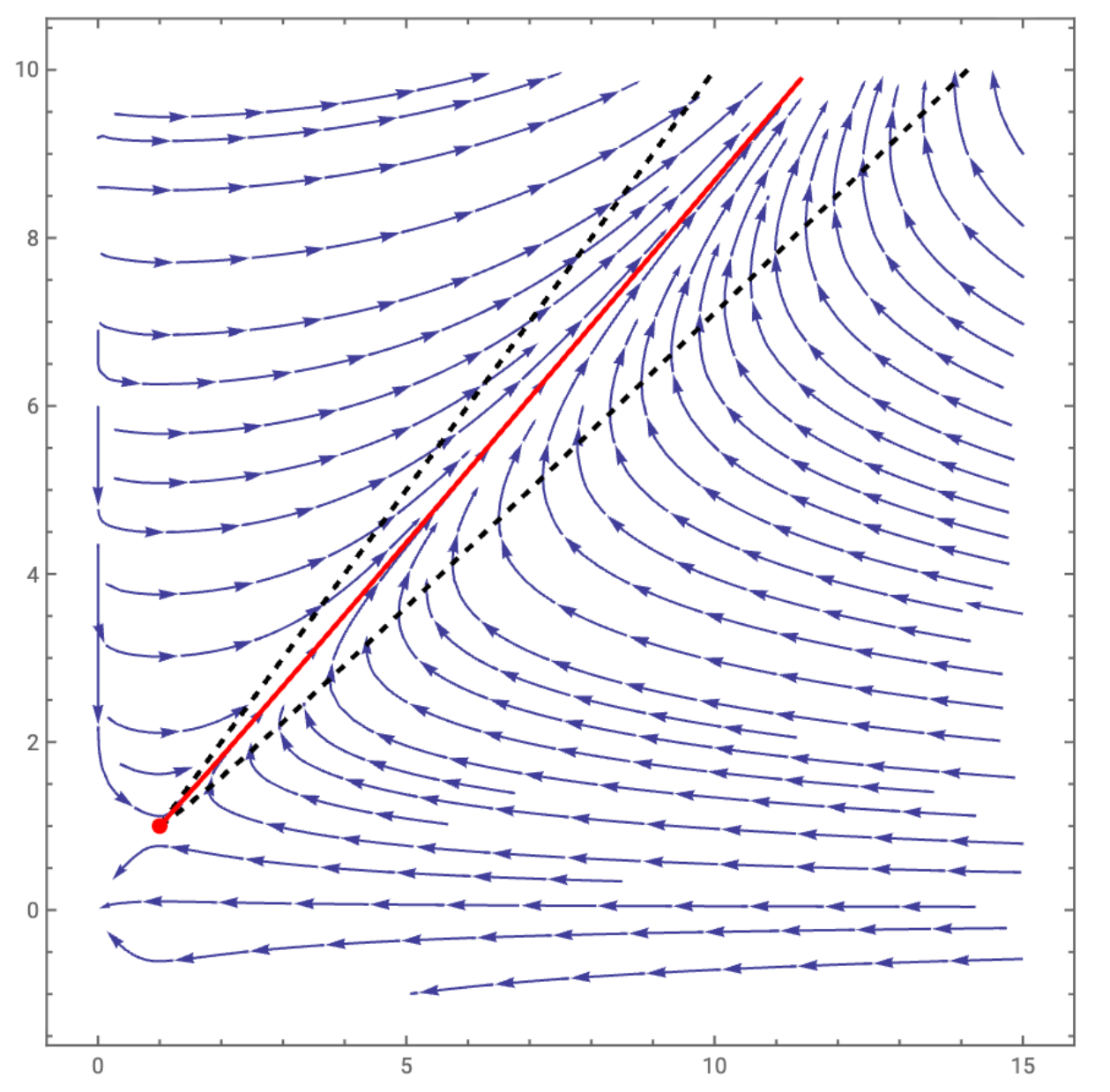}
\caption{{$(h,f)$-phase diagram of (\ref{eq:hfsystemALC}).}}
\label{fig:figure2}
\end{figure}
The phase diagram of \eqref{eq:hfsystemALC} is given in Figure \ref{fig:figure2}, where the red curve corresponds to a branch of the unstable manifold coming out of $(1,1)$. As in the previous case, we shall show again that this gives a complete strong $\G_2T$-structure on $S^3 \times S^3 \times \R$ though with a different asymptotic behaviour at one end.

\begin{proposition}\label{prop:ALCbranch}
The stationary point $(h,f)=(1,1)$ of \eqref{eq:hfsystemALC} is
hyperbolic, with both the stable and unstable manifolds 1-dimensional.
Moreover, one branch of the unstable manifold lies in the region
\[
 \mathcal R=
 \bigl\{(h,f):1<f<h,\quad h^2+1<2f^2\bigr\}
\]
and remains in $\mathcal R$ for all forward time.
\end{proposition}

\begin{proof}
The linearisation of \eqref{eq:hfsystemALC} at $(1,1)$ is
\[A=
 \begin{pmatrix}
 -2&4\\
 1&0
 \end{pmatrix},
\]
whose eigenvalues are
$ \lambda_\pm=-1\pm\sqrt5 \approx -3.24, 1.24.$ As before, from the stable/unstable manifold theorem we know that there are $1$-dimensional stable and unstable manifolds at $(1,1)$. We will consider the unstable branch (depicted in red in Figure \ref{fig:figure2}) which is tangent at $(1,1)$ to the line defined by $ h-1=\lambda_+(f-1)$ corresponding to the eigenvector 
$$
A\begin{pmatrix}h-1\\f-1\end{pmatrix}=\lambda_+\begin{pmatrix}h-1\\f-1\end{pmatrix}.
$$
From the linear approximation of the curve near $(1,1)$, it is not hard to see that this branch of the unstable manifold enters the region $\mathcal R$; this is also clear from Figure \ref{fig:figure2}. We now check that the branch cannot leave this region. On the boundary
$h=f>1$, we have
\[
h'=1-\frac1{h^2},
\qquad
f'=\frac12\left(1-\frac1{h^2}\right),
\]
and therefore
\[
(h-f)'=\frac12\left(1-\frac1{h^2}\right)>0.
\]
Thus, the branch cannot cross from $h>f$ to $h<f$. For the other boundary, set
$
G(h,f)=2f^2-h^2-1.
$
On $G=0$, we have $h'=0$, and hence
\[
G'
=4ff'-2hh'
=\frac{2(h^2-1)}{h}>0.
\]
Thus, the vector field also points into the region $G>0$. This concludes the proof.
\end{proof}

\begin{lemma}\label{lem:ALCasymptotics}
The branch in Proposition~\ref{prop:ALCbranch} is defined for all
$t\in\R$, and
$ h(t),f(t)\to+\infty$ as $t\to+\infty$. Moreover,
$ h(t)=\frac{2t}{3}+O(1)$ and $f(t)=\frac{t}{\sqrt3}+O(1)$
i.e. the corresponding metric is asymptotically locally conical.
\end{lemma}

\begin{proof}
The stable/unstable manifold theorem already gives exponential convergence to
$(h,f)=(1,1)$ as $t\to-\infty$. So it only remains to consider the positive end. Using \eqref{eq:hfsystemALC}, a direct calculation shows that inside the region $\mathcal R$ we have
\[
0<f'<\frac{h}{2f}<\frac1{\sqrt2},
\qquad
0<h'=2-\frac{h^2+1}{f^2}<2.
\]
Thus, $h$ and $f$ are strictly increasing and they cannot blow up in finite forward time. If they were bounded above then the 
trajectory would converge to a stationary point in the closure of
$\mathcal R$, but $(1,1)$ is the only such point and hence this cannot happen. The argument is analogous to that in Lemma \ref{lemma:ac} so we omit the details. Therefore, we conclude that $h(t),f(t)\to +\infty$. In particular, the solution is defined for all $t\in\R$, and the induced metric is complete.

We now determine the growth rate. Set
\[
q=\frac{h^2}{f^2},
\qquad
x=\log(f^2).
\]
A direct calculation using \eqref{eq:hfsystemALC} gives
\begin{equation}\label{eq:ALCq}
 \frac{dq}{dx}
 =
 \frac{q\bigl(f^2(4-3q)-1\bigr)}
 {qf^2-1}.
\end{equation}
Since we are in the region $\mathcal{R}$, the denominator is always positive. Moreover, \eqref{eq:ALCq} vanishes precisely on the barrier defined by
\[
q_*(f)=\frac43-\frac{1}{3f^2}.
\]
We note that the function $q_*$ is strictly increasing and converges to $\frac{4}{3}$ as
$f\to+\infty$. For a solution $q(f)<q_*(f)$, we have $\frac{dq}{dx}>0$ and if $q(f)>q_*(f)$ then $\frac{dq}{dx}<0$.

Near $(1,1)$, from the stable/unstable manifold theorem, the solution curve $(h(t),f(t))$ satisfies
\[
h(t)-1=\lambda_+(f(t)-1)+O((f(t)-1)^2).
\]
A direct calculation then shows
\[
q(t)
=1+2(\sqrt5-2)(f(t)-1)+O((f(t)-1)^2),
\]
and we also note that $2(\sqrt5-2)\approx 0.47$. 
On the other hand, we have the Taylor series expansion
\[
q_*(f)
=1+\frac{2}{3}(f-1)+O((f-1)^2).
\]
Thus, $q(t)<q_*$ near $(1,1)$ and hence is increasing. We claim that this holds for all $t$. At a hypothetical first point of
contact $x_0$ where $q(t)$ meets $q_*$, we would have $\frac{dq}{dx}\big|_{x_0}=0$, and also since $q(t)<q_*(f)$ for $x<x_0$ it follows that $\frac{d}{dx}(q-q_*)\big|_{x_0}\geq 0$. However, we saw that $q_*(f)$ is an increasing function and as such $\frac{d}{dx}q_*>0$; hence this gives a contradiction. So we conclude that 
the solution cannot cross the $q_*$ from below.
Therefore,
\[
1<q(t)<q_*(f)<\frac43,
\]
and $q(t)$ is increasing.

Let $L\leq\frac43$ denote the limit of $q(t)$. If $L<\frac43$, then the
right-hand side of \eqref{eq:ALCq} tends to the positive constant
$4-3L$, which is impossible since $q$ converges. Therefore,
\begin{equation}\label{eq:ALCratio}
 \frac hf\longrightarrow\frac2{\sqrt3}.
\end{equation}
Substituting in \eqref{eq:hfsystemALC} gives
\[
h'(t)\longrightarrow\frac23,
\qquad
f'(t)\longrightarrow\frac1{\sqrt3},
\qquad
h(t)\sim\frac{2t}{3},
\qquad
f(t)\sim\frac{t}{\sqrt3}.
\]
To obtain the decay rate, similarly as before we first set
$ u=q-\frac43. $
Equation \eqref{eq:ALCq} then becomes
\[
u_x+3u
=
-\frac{\left(\frac43+4u\right)f^{-2}}
{\frac43+u-f^{-2}}.
\]
Since $u\to0$ and $f^{-2}=e^{-x}\to0$, it follows that
\[
u_x+3u=O(e^{-x}).
\]
Thus,
\[
\frac{d}{dx}(e^{3x}u)=O(e^{2x}),
\]
and integration gives
\[
u=O(e^{-x})=O(f^{-2})=O(t^{-2}).
\]
Therefore,
\[
\frac{h^2}{f^2}
=\frac43+O(t^{-2}).
\]
Substituting  in \eqref{eq:hfsystemALC}, we obtain
\[
h'
=2-\frac{h^2}{f^2}-\frac1{f^2}
=\frac23+O(t^{-2})\qquad 
f'
=\frac{\frac{h^2}{f^2}-\frac1{f^2}}
{2\frac hf}
=\frac1{\sqrt3}+O(t^{-2}).
\]
Hence
\[
h(t)=\frac{2t}{3}+O(1),
\qquad
f(t)=\frac{t}{\sqrt3}+O(1).
\]
Thus, the corresponding ALC model is given by
\begin{equation}\label{eq:ALCmodel}
 g_{\mathrm{ALC}}
 =
 dt^2+4(\eta^2)^2
 +\frac{16}{9}t^2(\eta^1)^2
 +\frac43t^2
 \bigl((\eta^3)^2+(\eta^4)^2+(\eta^5)^2+(\eta^6)^2\bigr).
\end{equation}
Here $\eta^2$ corresponds to the connection $1$-form of the constant-length circle
direction, while the remaining part corresponds to the Calabi-Yau cone over
$T^{1,1}\cong S^2\times S^3$. From the above estimates we have
$
|g_\varphi-g_{\mathrm{ALC}}|_{g_{\mathrm{ALC}}}
=O(t^{-1})
$ and the corresponding higher derivatives estimates follow by
differentiating the ODE. Thus, the positive end is ALC with rate
$\nu=-1$.
\end{proof}

The above can be summarised as follows:
\begin{theorem}\label{thm:mainALC}
There exists a complete strong $\G_2T$-structure on
$\R\times S^3\times S^3$ whose associated metric is asymptotically
cylindrical at one end, asymptotic
to the product of the Euclidean metric on $\R$ with a bi-invariant
metric on $S^3\times S^3$, and is asymptotically locally conical at the other. The ALC model metric is given by \eqref{eq:ALCmodel}\footnote{With the normalisation $A_1=a$, the constant circle direction has length $2a$. Without loss of generality, we set $a=1$ in the above.}. The characteristic
torsion is independent of $t$ and is given by
\[
T_\varphi=-4(e^{123}+e^{456}).
\]
\end{theorem}

One can again check directly, using say Mathematica, that the strong $\G_2T$ metric obtained in the above theorem has full holonomy group i.e. $\mathrm{dim}(\mathfrak{hol}(\nabla^{LC}))=21$. Moreover, as with the previous example, one can again verify that the Lee form $\theta=-\frac{2}{hf^2}dt$ is not Killing and hence the characteristic Ricci tensor is non vanishing.

\section{The quotient \texorpdfstring{$\SU(3)$}{}-structure}\label{sec: quotient}

We now investigate the $\SU(3)$-structure obtained by taking a quotient of the examples constructed in the previous section. This is motivated by our results in \cite{FinoFowdar2026}*{\S 5} and the observation that each of the complete strong $\G_2 T$-structures given in Theorem \ref{thm:mainsu(2)3} and \ref{thm:mainALC} is invariant under an $S^1$ action. This circle action is generated by the left invariant vector field $X:=e_1+e_4$ satisfying $\eta^2(X)=1$ and $\eta^i(X)=0$ for $i\neq 2$; unlike the reduction considered in Section \S 5, this action is \textit{not} generated by the Lee vector field.
It turns out that the quotient $\SU(3)$-structure has a surprisingly rich structure; though they are very different in each case. For the AC example, we shall show that the almost Hermitian structure is conformally almost K\"ahler and the almost complex structure is non-integrable, while for the ALC example we get a complex manifold
which is conformally K\"ahler. In the latter case, the Riemannian 
connection $1$-form for the circle action is Hermitian-Yang-Mills and  the Hermitian structure is Calabi-Yau with torsion i.e. its Bismut Ricci form vanishes. We now describe these results in more detail. 

The $1$-form $\eta:=\eta^2$ defines a connection for the principal circle bundle generated by the vector field $X$
\[
 S^1\hookrightarrow \mathbb R\times S^3\times S^3
 \longrightarrow P^6:=\mathbb R\times (S^3\times S^3)/S^1
 \cong\mathbb R\times S^2\times S^3,
\]
with curvature form given by
\begin{equation}\label{eq:quotient-curvature}
 F_\eta=d\eta=-2(\eta^{35}+\eta^{46}).
\end{equation}
It is easy to check that $\mathcal{L}_X\varphi=0$ and hence the quotient
inherits an $\SU(3)$-structure $(\omega,\Omega^\pm)$. We denote by
\[
 u:=g_{\vp}(X,X)^{-1/2}=\frac{1}{2A_1}
\]
the (inverse) length of the circle fibres; this was denoted by $t$ in \cite{FinoFowdar2026}*{\S 5} but to avoid confusion with the radial coordinate, we shall use $u$ instead. Following the construction in \cite{FinoFowdar2026}, we set
\begin{equation}\label{eq:quotient-reduction}
 \vp=u^{-3}(\eta\wedge\omega+\Omega^+),\qquad
 *_{\vp}\vp=u^{-4}\left(\frac12\omega^2-\eta\wedge\Omega^-\right),
 \qquad g_{\vp}=u^{-2}(\eta\otimes\eta+g_\omega).
\end{equation}
For the $\G_2$-structure defined by \eqref{eq: vpansatz}, the quotient $\SU(3)$-structure can be described explicitly as follows:
\begin{proposition}\label{prop:quotient-SU3}
For the general four-function ansatz \eqref{eq: vpansatz}, the induced
$\SU(3)$-structure is given by
\begin{align}
 \omega={}&-\frac{B_1}{2A_1^2}\,\eta^1\wedge dt
   -\frac{A_2B_2}{A_1^2}(\eta^{36}+\eta^{45}),
   \label{eq:quotient-omega}\\
 \Omega^+={}&\frac{B_1}{A_1^3}
       (B_2^2\eta^{135}-A_2^2\eta^{146})
       +\frac{A_2B_2}{2A_1^3}
       (\eta^{34}+\eta^{56})\wedge dt,
       \label{eq:quotient-Omegaplus}\\
 \Omega^-={}&\frac{B_1A_2B_2}{A_1^3}
       (\eta^{134}+\eta^{156})
       +\frac{1}{2A_1^3}
       (A_2^2\eta^{46}-B_2^2\eta^{35})\wedge dt.
       \label{eq:quotient-Omegaminus}
\end{align}
In particular, the associated metric is
\begin{equation}\label{eq:quotient-metric}
 \begin{aligned}
 g_\omega={}\frac{dt^2}{4A_1^2}
 +\frac{B_1^2}{A_1^2}(\eta^1)^2
 +\frac{B_2^2}{A_1^2}\big((\eta^3)^2+(\eta^5)^2\big)
 +\frac{A_2^2}{A_1^2}\big((\eta^4)^2+(\eta^6)^2\big).
 \end{aligned}
\end{equation}
\end{proposition}

The next observation distinguishes the ALC family from the other
solutions of the four-function system. As before, we shall always assume that all four
coefficient functions $A_i,B_i$ are positive.

\begin{proposition}\label{prop:quotient-integrability}
Let $\vp$ be a strong $\G_2T$-structure satisfying
\eqref{eq:A1system}-\eqref{eq:B2system}. The almost complex structure on
the $S^1$ quotient is integrable if and only if
\[
 \kappa=-A_1^2,
\]
which in particular forces $A_1$ to be constant. 
Moreover, the curvature $2$-form $F_\eta$ is of type $(1,1)$ if and only if $A_2=B_2$, and satisfies $\iota_{\omega}F_\eta=0$.
\end{proposition}
\begin{proof}
Using the orthonormal co-framing: 
\[
 \epsilon^1=\frac{dt}{2A_1},\ \
 \epsilon^2=\frac{B_1}{A_1}\eta^1,\ \
 \epsilon^3=\frac{B_2}{A_1}\eta^3,\ \
 \epsilon^4=-\frac{A_2}{A_1}\eta^6,\ \
 \epsilon^5=\frac{A_2}{A_1}\eta^4,\ \
 \epsilon^6=-\frac{B_2}{A_1}\eta^5,
\]
on $P^6$, we can express the $\SU(3)$-structure in the standard form:
\[
 \omega=\epsilon^{12}+\epsilon^{34}+\epsilon^{56},\qquad
 \Omega=(\epsilon^1+i\epsilon^2)
             \wedge(\epsilon^3+i\epsilon^4)
             \wedge(\epsilon^5+i\epsilon^6).
\]
Denoting the complex $(1,0)$ forms by $\zeta^1=\epsilon^1+i\epsilon^2$, $ \zeta^2=\epsilon^3+i\epsilon^4$, 
$ \zeta^3=\epsilon^5+i\epsilon^6$, 
a direct calculation using the structure equations and the ODE system gives:
\begin{align*}
 (d\zeta^1)^{0,2}&=0,\\
 (d\zeta^2)^{0,2}&=\frac{(A_1^2+\kappa)(A_2^2+B_2^2)}
 {4A_2^2B_2^2}\,\overline\zeta^1\wedge\overline\zeta^2,\\
 (d\zeta^3)^{0,2}&=-\frac{(A_1^2+\kappa)(A_2^2+B_2^2)}
 {4A_2^2B_2^2}\,\overline\zeta^1\wedge\overline\zeta^3.
\end{align*}
Thus, the complex structure is integrable if and only if
$\kappa=-A_1^2$. Since $\kappa$ is constant, this forces $A_1$ to be
constant; then \eqref{eq:A1system} implies $A_2^2=B_2^2$, and
hence $A_2=B_2$. From \eqref{eq:quotient-curvature} we also have
\[
 F_\eta=2A_1^2\left(\frac{1}{B_2^2}\epsilon^{36}
                      -\frac{1}{A_2^2}\epsilon^{45}\right).
\]
and thus,  
\[
 (F_\eta)^{(2,0)+(0,2)}=
 A_1^2\left(\frac{1}{B_2^2}-\frac{1}{A_2^2}\right)
 (\epsilon^{36}+\epsilon^{45}).
\]
This shows that $F_\eta$ is of type $(1,1)$ precisely when $A_2=B_2$ and it is orthogonal with respect to $\omega$. This concludes the proof.
\end{proof}

\subsection{The quotient of the complete AC example}
We now focus on the $\SU(2)^3$-invariant strong $\G_2 T$-structure. In the notation of \S7.2.1, the quotient $\SU(3)$-structure for the example given in  Theorem~\ref{thm:mainsu(2)3} can expressed as:
\begin{align}\label{eq:quotient-AC-metric}
 g_\omega={}&\frac{dt^2}{f^2}
       +\frac{h^2}{f^2}\big((\eta^1)^2+(\eta^3)^2+(\eta^5)^2\big)
       +(\eta^4)^2+(\eta^6)^2,\\
 \omega={}&-\frac{h}{f^2}\eta^1\wedge dt
       -\frac hf(\eta^{36}+\eta^{45}),
       \label{eq:quotient-AC-omega}\\
 \Omega^+={}&\frac{h^3}{f^3}\eta^{135}
       -\frac hf\eta^{146}
       +\frac h{f^2}(\eta^{34}+\eta^{56})\wedge dt,
       \label{eq:quotient-AC-Omegaplus}\\
 \Omega^-={}&\frac{h^2}{f^2}(\eta^{134}+\eta^{156})
       +\left(\frac1f\eta^{46}-\frac{h^2}{f^3}\eta^{35}\right)
       \wedge dt.
       \label{eq:quotient-AC-Omegaminus}
\end{align}
A direct calculation shows that the only non-zero $\SU(3)$ torsion forms are:
\begin{align}
 \nu_1&=\frac{3f^4-3f^2h^2-f^2-h^2}{f^3h^2}\,dt,
       \label{eq:quotient-AC-nu1}\\
 \pi_1&=\frac{2(2f^4-2f^2h^2-f^2-h^2)}{f^3h^2}\,dt,
       \label{eq:quotient-AC-pi1}\\
 \pi_2&=\frac{(f^2-1)(f^2+h^2)}{f^2h^2}
                   (\epsilon^{36}-\epsilon^{45}),
       \label{eq:quotient-AC-pi2}
\end{align}
where $\epsilon^i$ denotes the orthonormal quotient co-framing
introduced in Proposition~\ref{prop:quotient-integrability}. 
Thus, since $\pi_2\neq 0$ the induced almost complex structure is not integrable and its
Nijenhuis tensor is not totally skew-symmetric, see \cite{FinoFowdar2026}*{Proposition 4.1}. On the other hand, it is clear that 
$\nu_1$ is exact and hence
$$d\omega=\nu_1\wedge\omega,$$ 
so the quotient
is globally conformally almost K\"ahler. It is also clear that the curvature form
\[
 F_\eta=2\left(\frac{f^2}{h^2}\epsilon^{36}-\epsilon^{45}\right),
\]
is not of type $(1,1)$ though $F_\eta\wedge\omega^2=0$. 

We now look at the two asymptotics of the quotient structure. We begin by introducing the geodesic coordinate:
\[
 s(t)=\int_{t_0}^t\frac{d\xi}{f(\xi)},
\]
so that the metric \eqref{eq:quotient-AC-metric} can be expressed as
\[
 g_\omega=ds^2
       +\frac{h^2}{f^2}\big((\eta^1)^2+(\eta^3)^2+(\eta^5)^2\big)
       +(\eta^4)^2+(\eta^6)^2.
\]
From our previous analysis, we have 
$h,f\to1$ as $t\to-\infty$, and as $t\to+\infty$, we have
\[
 f(t)\sim\frac23t,\qquad
 \frac{h(t)}{f(t)}\longrightarrow\sqrt3,
 \qquad s(t)=\frac32\log t+O(1).
\]
Thus, the quotient metric $g_{\omega}$ is complete and has two asymptotically cylindrical ends. As $t\to -\infty$, we have
\[
g_{\omega}\sim ds^2 + (\eta^1)^2+(\eta^3)^2+(\eta^4)^2
                +(\eta^5)^2+(\eta^6)^2
\]
and as $t\to +\infty$, we have
\[
 g_{\omega}\sim ds^2+ 3\big((\eta^1)^2+(\eta^3)^2+(\eta^5)^2\big)
                +(\eta^4)^2+(\eta^6)^2.
\] 
In particular, we see that although the metric on $M^7=\R\times S^3 \times S^3$ has an AC end, its circle quotient is asymptotically cylindrical at that end.

\subsection{The quotient of the complete ALC example}

We now consider the example obtained in Theorem~\ref{thm:mainALC}. Recall that in this case we had normalised $A_1=-\kappa=1$ so that the $S^1$-fibres have constant length $2$. In the notation of \S7.2.2, the quotient
$\SU(3)$-structure is given by
\begin{align}
 g_\omega &=\frac14dt^2+h^2(\eta^1)^2
       +f^2\big((\eta^3)^2+(\eta^4)^2
                    +(\eta^5)^2+(\eta^6)^2\big),\label{eq:quotient-ALC-metric}\\
 \omega&=-\frac h2\eta^1\wedge dt
             -f^2(\eta^{36}+\eta^{45}),
             \label{eq:quotient-ALC-omega}\\
 \Omega^+&=hf^2(\eta^{135}-\eta^{146})
          +\frac{f^2}{2}(\eta^{34}+\eta^{56})\wedge dt,
          \label{eq:quotient-ALC-Omegaplus}\\
 \Omega^-&=hf^2(\eta^{134}+\eta^{156})
          +\frac{f^2}{2}(\eta^{46}-\eta^{35})\wedge dt.
          \label{eq:quotient-ALC-Omegaminus}
\end{align}
As $t\to-\infty$, it is easy to see that metric is asymptotic to the homogeneous metric:
$$g_{\omega}\sim \frac{1}{4}dt^2+(\eta^1)^2+(\eta^3)^2+(\eta^4)^2+(\eta^5)^2+(\eta^6)^2.$$
On the other hand, as $t\to +\infty$, from the previous analysis of the asymptotics of $h$ and $f$ we have:
\[
 g_\omega\sim \frac{1}{4}dt^2+\frac{1}{4}t^2\left[\frac{16}{9}(\eta^1)^2
  +\frac43\big((\eta^3)^2+(\eta^4)^2
                    +(\eta^5)^2+(\eta^6)^2\big)\right],
\]
which correpsonds to the conical  Calabi-Yau metric over $T^{1,1}=S^3\times S^2$. This is precisely the \textit{conifold metric}. 

The torsion forms of the induced $\SU(3)$-structure on
$P^6\cong\mathbb R\times S^2\times S^3$ are given by:
\begin{equation}\label{eq:quotient-differentials}
 d\omega=\alpha\wedge\omega,\qquad
 d\Omega^+=2\alpha\wedge\Omega^+,\qquad
 d\Omega^-=2\alpha\wedge\Omega^-
\end{equation}
where  $\alpha=-{h^{-1}f^{-2}}dt$ is an exact $1$-form. Thus, it is clear that $\pi_0=\sigma_0=\pi_2=\sigma_2=0$ i.e. its almost complex structure is integrable. Moreover, since $\alpha$ is exact, the underlying 
Hermitian metric is globally conformally K\"ahler. 

As the underlying $\SU(3)$-structure is Hermitian, it admits a Bismut connection. Thus, it is natural to ask if it satisfies any interesting geometric properties.
\begin{proposition}\label{prop:quotient-bianchi}
For the quotient of the complete ALC strong $\G_2T$-structure,
the circle connection $\eta=\eta^2$ is Hermitian-Yang-Mills, i.e.
\begin{equation}\label{eq:quotient-HYM}
 F_\eta\wedge\Omega^+=F_\eta\wedge\Omega^-=0,
 \qquad F_\eta\wedge\omega^2=0.
\end{equation}
The Bismut torsion form of $(P,\omega,\Omega^{\pm})$ is given by
\begin{equation}\label{eq:quotient-Bismut-torsion}
 T_\omega=-2\eta^1\wedge(\eta^{36}+\eta^{45})
           =\eta^1\wedge d\eta^1,
\end{equation}
and satisfies the Bianchi identity:
\begin{equation}\label{eq:quotient-Bianchi}
 dT_\omega=-F_\eta\wedge F_\eta.
\end{equation}
Moreover, the Bismut-Ricci form  $\rho^B$ vanishes, i.e. the underlying structure is Calabi-Yau with torsion, but it is not pluriclosed.
\end{proposition}
\begin{proof}
    This follows from a straightforward calculation. To see that $\rho^B$ vanishes, it suffices to note that the $\SU(3)$ torsion forms satisfy $\pi_1=2\nu_1=2\alpha$ and hence the result follows from \cite{FinoFowdar2026}*{Corollary 4.8}.
\end{proof}
\begin{remark}
The relation between the $7$-dimensional characteristic torsion and
the Bismut torsion of the quotient is given by
\begin{equation}\label{eq:quotient-torsion-relation}
 T_{\vp}=4(T_\omega+\eta\wedge F_\eta).
\end{equation}
We emphasise that the circle action considered here is not generated by the $\G_2$ Lee vector field. In fact, for the complete ALC example
we have $\theta=4\tau_1=2\alpha$, which is radial and non-parallel.
Consequently, this quotient is different from the characteristic
Ricci-flat reduction considered in \cite{FinoFowdar2026}*{\S 6}.
\end{remark}

The key differences between the two $\SU(3)$ quotients can be summarised as follows:
\begin{table}[htbp]
\centering
\renewcommand{\arraystretch}{1.25}
\begin{tabularx}{\textwidth}{@{}lXX@{}}
\toprule
\textbf{Property}
& \textbf{AC quotient}
& \textbf{ALC quotient} \\
\midrule
Underlying manifold
& $\mathbb{R}\times S^2\times S^3$
& $\mathbb{R}\times S^2\times S^3$ \\

Complete metric
& Yes
& Yes \\

Almost complex structure
& Non-integrable and $N_J \notin \Lambda^3$
& Integrable \\

Conformal geometry
& Conformally almost K\"ahler
& Conformally K\"ahler \\

Negative end
& Cylindrical
& Cylindrical \\

Positive end
& Cylindrical
& Conical Calabi-Yau \\
\bottomrule
\end{tabularx}
\caption{Comparing the induced $\mathrm{SU}(3)$-structures
on the quotients of the AC and ALC strong
$\mathrm{G}_2T$-structures.}
\label{tab:quotient-comparison}
\end{table}
Note that since the Nijenhuis tensor $N_J$ is not totally skew-symmetric for the AC quotient and hence the underlying almost Hermitian structure on $P$ does \textit{not} admit a Hermitian connection with skew-symmetric torsion \cite{FriedrichIvanov2002}.

\bibliographystyle{plain}
\bibliography{biblioG}

\noindent
A. Fino \\
Dipartimento di Matematica ``G. Peano'', Universit\`{a} degli studi di Torino \\
Via Carlo Alberto 10, Torino (Italy) \\
\& Department of Mathematics and Statistics, Florida International University\\
Miami, FL 33199 (United States)\\
Email: \texttt{annamaria.fino@unito.it, afino@fiu.edu}

\vskip0.8em
\noindent
U. Fowdar \\ 
Institute of Mathematics, University of Warsaw\\ 
Banacha 2, 02-097 Warszawa (Poland)\\
Email: \texttt{u.fowdar@uw.edu.pl} 

\end{document}